%% file: ARTICLE.tex
\documentclass[preprint,10pt]{elsarticle}

\usepackage{amssymb}

\usepackage[utf8]{inputenc}
\usepackage{amsmath}
\usepackage{amsfonts}
\usepackage{enumitem}  
\usepackage[ruled,vlined,linesnumbered,resetcount]{algorithm2e}
\usepackage{tikz}
\usetikzlibrary{positioning}
\usepackage{caption}
\usepackage{subcaption}

\usepackage{xcolor}
\usepackage{graphicx}
\usepackage[normalem]{ulem}

\usepackage{bm}

\usetikzlibrary{arrows.meta}

\usepackage{minibox}
\usepackage[normalem]{ulem}

\journal{International Journal for Numerical Methods in Fluids}

\input{./NEW_COMMANDS.tex}

\usepackage[a4paper, total={5.5in, 9.in}]{geometry}

\usepackage{floatrow}
\newfloatcommand{capbtabbox}{table}[][\FBwidth]

\SetKwInput{KwOffline}{Offline}
\SetKwInput{KwOnline}{Online}
\usepackage[export]{adjustbox}
\begin{document}

\begin{frontmatter}



\author{M. Oulghelou\fnref{label1}}
\fntext[label1]{mourad.oulghelou@univ-lr.fr}
\author{C. Allery\fnref{label2}}
\fntext[label2]{cyrille.allery@univ-lr.fr}

\address{LaSIE, UMR-7356-CNRS, Universit\'e de La Rochelle P\^ole Science et Technologie,
Avenue Michel Cr\'epeau, 17042 La Rochelle Cedex 1, France.}



\title{Hyper bi-calibrated interpolation on the Grassmann manifold for near real time flow control using genetic algorithm}


\begin{abstract}
\input{./Tex_Files/abstract.tex}
\end{abstract}


\end{frontmatter}
\newpage
\section{Introduction}
\input{./Tex_Files/Introduction}

\section{Bi-CITSGM method}
\input{./Tex_Files/intro_Review_Bi-ITSGM}
\input{./Tex_Files/Subspaces_interpolation}
\subsection{Bi-CITSGM method}
\input{./Tex_Files/Enhanced_Bi-CITSGM}

\section{Hyper Bi-CITSGM}
In this section, we propose the Hyper Bi-CITSGM approach designed to accelerate the CPU time of the usual Bi-CITSGM approach suggested in \cite{OulghelouBiCITSGM2019Arxiv}. The description of this approach is given in th following. 
\subsection{Description of the approach}
\input{./Tex_Files/Hyper_Bi-CITSGM}
\subsection{Computational complexity}
\input{./Tex_Files/Comp_complexity_Hyper_Bi-CITSGM.tex}

\section{Near real time flow control using a reduced genetic algorithm}
	\subsection{Problem setting}
		\input{./Tex_Files/ctrl_prob.tex}	
	\subsection{Canonical Genetic algorithm}	
		\input{./Tex_Files/GA_overview.tex}

	\subsection{Reduced Genetic Algorithm}	
		\input{./Tex_Files/used_GA.tex}	

\section{Numerical applications}	
	\subsection{Control of flow past a cylinder}
		\input{./Tex_Files/GA_Cylinder.tex}
	\subsection{Control of flow in a lid driven cavity}
			\input{./Tex_Files/GA_Cavity.tex}
\section{Conclusions}
	\input{./Tex_Files/Conclusion.tex}
\section*{Acknowledgement}
The authors would like to thank the Nouvelle-Aquitaine region and the European union for their financial support. This material is based upon work financially supported by CPER BATIMENT DURABLE - Axe 3 "Qualité des Environnement Intérieurs (QEI)", Convention number : P-2017-BAFE-102.

%
%

\section*{\large References}
\bibliographystyle{ieeetr}
\bibliography{./ARTICLE}
\end{document}

%% file: NEW_COMMANDS.tex
\newcommand{\myfrac}[2]{\displaystyle{\frac{#1}{#2}}}
\newcommand{\somme}[3]{\displaystyle{\overset{#3}{\underset{#1=#2}{\sum}}}}

\newcommand{\vect}[1]{ {\bm{#1}}}

\newcommand{\ctrlvar}{\gamma}

\newcommand{\statevecU}{\vect{u}}
\newcommand{\statevecUCible}{\hat{\vect{u}}}

\newcommand{\InterpSbase}{{\Phi_{_{\untrainedparam}}}}
\newcommand{\InterpTbase}{{\Lambda_{_{\untrainedparam}}}}

\newcommand{\Sbase}[1]{{\Phi_{_{\trainingparam{#1}}}}}

\newcommand{\Xbasis}[1]{{\Phi_{_{#1}}}}
\newcommand{\XbasisEnd}[1]{{\Psi_{_{#1}}}}

\newcommand{\Tbase}[1]{{\Lambda_{_{\trainingparam{#1}}}}}
\newcommand{\trainedSingVal}[1]{{\Sigma_{_{\trainingparam{#1}}}}}
\newcommand{\InterpSingVal}{{\Sigma_{_{\untrainedparam}}}}

\newcommand{\nbrModes}{q}

\newcommand{\nbrParam}{N_{p}}
\newcommand{\nbrSnap}{N_{s}}
\newcommand{\nbrBases}{N_{p}}
\newcommand{\dimVecorsSnap}{N_x}
\newcommand{\dimParam}{p}

\newcommand{\trainedSnapMat}[1]{\bm{S}_{_{\trainingparam{#1}}}}
\newcommand{\trainedredSnapMat}[1]{S^{(r)}_{_{ \trainingparam{#1}}}}

\newcommand{\InterpSnapMat}{\bm{S}_{_{\untrainedparam}}}
\newcommand{\InterpRedSnapMat}{S^{(r)}_{_{ \untrainedparam}}}

\newcommand{\mcal}[1]{\mathcal{#1}}
\newcommand{\mbb}[1]{\mathbb{#1}}

\newcommand{\txt}[1]{\textnormal{#1}}

\newcommand{\constraint}[0]{\mathcal{N}}

\newcommand{\functional}[0]{\mathcal{J}}

\newcommand{\GrassmanifoldSpace}[1]{\mathcal{G}(\nbrModes, N_{#1})}

\newcommand{\TangsubSpaceGrass}[2]{\mcal{T}_{_{[{\Mat{#1}}_{#2}]}}\GrassmanifoldSpace{}}

\newcommand{\OrthGroupe}[1]{\bm{\mcal{O}}(#1)}

\newcommand{\Exp}[1]{\bm{\txt{Exp}}_{_{[#1]}}}
\newcommand{\Log}[1]{\bm{\txt{Log}}_{_{[#1]}}}

\newcommand{\Mat}[1]{{#1}}
\newcommand{\MatCalib}[1]{{#1}}

\newcommand{\distGrass}[2]{\txt{\textit{dist}}_{\mcal{G}} (#1,#2)}

\newcommand{\trainingparam}[1]{\ctrlvar_{#1}}
\newcommand{\untrainedparam}{\widetilde{\ctrlvar}}

\newcommand{\InitVel}[1]{\mcal{X}_{_{#1}}}

\newcommand{\InterpGrassVel}{\mcal{X}_{_{\untrainedparam}}}

\newcommand{\CalibMatrix}[2]{\Mat{#1}_{_{#2}}}

%% file: Tex_Files/abstract.tex
Most genetic algorithms (GAs) used in literature to solve control problems are time consuming and involve important storage memory requirements. In fact, the search in GAs is iteratively performed on a population of chromosomes (control parameters). As a result, the cost functional needs to be evaluated through solving the high fidelity model or by performing the experimental protocol for each chromosome and for many generations. To overcome this issue, a non intrusive reduced real coded genetic algorithm (RGA) for near real time optimal control is designed. 
This algorithm uses precalculated parametrized solution snapshots stored in the POD (Proper Orthogonal Decomposition) reduced form, to predict the solution snapshots for chromosomes over generations. The method used for this purpose is a hyper reduced version of the Bi-CITSGM method (Bi Calibrated Interpolation on the Tangent Space of the Grassmann Manifold) designed specially for non linear parametrized solution snapshots interpolation \cite{OulghelouBiCITSGM2019Arxiv}. This hyper reduced approach referred to as Hyper Bi-CITSGM, is proposed in such a way to accelerate the usual Bi-CITSGM process by bringing this last to a significantly low dimension.
Thus, the whole optimization process by RGA can be performed in near real time. The potential of RGA in terms of accuracy and CPU time is demonstrated on control problems of the flow past a cylinder and flow in a lid driven cavity when the Reynolds number value varies.

%% file: Tex_Files/Introduction.tex
The most used optimization techniques for flow optimal control are of gradient descent type \citep{Gunzburger2000, Snyman2005}. Considerable research efforts were conducted in the development of techniques for evaluating the sensitivity of the cost function with respect to optimization parameters. In particular, Lagrange theory played a major role in this subject and was widely used for constrained optimization problems \cite{Desai1994}. Using the Lagrangian approach, the Fréchet derivative of the cost functional can directly be determined via the solution of an auxiliary adjoint partial differential equation. Precisely, the solutions of two PDEs (state and adjoint equations) are needed each time the descent direction has to be updated. 
Eventhough their practical effectiveness, gradient-based strategies are susceptible to generate spurious local minima, which may inhibit their capabilities in some flow control applications \cite{Zingg2008}.
To circumvent this limitation, some authors proposed the employment of more powerful optimization strategies, such as Genetic Algorithms (GAs). This algorithm developed by John Holland \cite{Holland1975}, is a stochastic optimization approach known by its ability to perform global optimization \cite{Goldberg1989}. 
It is a form of evolutionary search that makes use of genetic operators that mimic the Darwinian concept of natural selection and evolution. These are typically selection, crossover and mutation. The flow optimal control solutions are coded in such a way that they can be thought of as forms of genetic material (DNA). A population of solutions is generated randomly and the fitness of each chromosome is assessed in such a way fit chromosomes possess greater chance of reproducing and thus promoting their fitter characteristics through to the subsequent generations. Crossover takes place by swapping parts of the DNA from two chromosomes and mutation by randomly alter genes of a chromosome.
Given that fitness of chromosomes is directly evaluated in terms of optimization parameters without the need of local derivatives information, GAs can deal with non-smooth, non-continuous and non-differentiable functions which are actually encountered in practical optimization situations.
Many attempts have been made with the purpose of using GAs for flow optimization. For instance, GAs were applied by Sengupta et al. \cite{SENGUPTA2007} to optimally control incompressible viscous flow past a circular cylinder for drag minimization by rotary oscillation; by Hacio\u{g}lu et al. \cite{Hacioglu2005} and Shahrokhi et al. \cite{SHAHROKHI2007} for airfoil shape optimization; and by Dar\`oczy et al. \cite{Daroczy2018} for the optimization of an H-Darrieus wind turbine. Despite their robustness, a serious weakness of GAs is their substantial lack of computational efficiency \cite{Vicini1999, TerryAerodynamicSO}. 
Concretely, the number of the cost function evaluations required by GAs exceeds in general the number required by a gradient-based optimization \cite{TerryAerodynamicSO, Shigeru1997}. For that reason, a great effort is made by scientists to accelerate GA-based optimization methods. Various solutions to accelerate and improve the performance of GAs were then suggested. For instance, by using improved genetic operators \cite{Hacioglu2002, Hacioglu2003}, by using of multiprocessing \cite{Doorly1998, Jones2000}, or by hybridization of GAs with a descent optimization method \cite{Vicini1999}. 
Here, we rather focus on the fitness evaluation stage where the costly high fidelity model is needed to be solved for each chromosome in the population. A possible way to drastically reduce the computational cost of the fitness evaluation is by using reduced order models (ROMs). 
\\
The most popular method in model order reduction is the Proper Orthogonal Decomposition (POD) \cite{lumley1967}. Starting from a set of solution snapshots, The POD method can generate an optimal basis in the sens that only few modes will be sufficient to reproduce the dynamics of the problem. 
%
%
%
However, a weak point of this method is its sensitivity to parameters changes. Thus, a POD reduced order model generated for a set of trained parameters cannot be expected to approximate well the dynamics for a new untrained parameter.
A possible way to circumvent this issue is by using the strategy proposed by Amsallem and Farhat for POD bases adaptation named here ITSGM (Interpolation on the Tangent Space of the Grassmann Manifold) \cite{Amsallem}. This method is a generalization of the Subspace Angle Interpolation (SAI) \cite{Lieu2004, LIEU20065730} relying on the concepts of principal angles between two subspaces and principal vectors for a pair of subspaces \cite{Bjorck71numericalmethods}.
In the context of optimal control using reduced order models, the ITSGM was successfully applied in the control of the non linear heat and Burgers equations \cite {OULGHELOUAMC2018}. 
It is worth mentioning that in this classical reduced control approach, the reduced order model describing the temporal dynamics is constructed via the Galerkin projection. This approach is thus considered intrusive, given the fact that it requires the access to the high fidelity model. Contrarily to this last, a different situation of optimal control where no prior required knowledge of the governing equations is studied in this article. This situation suggests to optimally control the physical problem by using only available solutions data of the corresponding high fidelity model in a bunch of trained control parameters.
For this purpose, we use a non intrusive reduced order model approach based on the ITSGM method.
This approach named Bi-CITSGM (Bi-Calibrated Interpolation on the Tangent Space of the Grassmann Manifold) has been recently developed and successfully applied to predict in real time the flow past a cylinder when the Reynolds number value varies \cite{OulghelouBiCITSGM2019Arxiv}. By using a set of precalulated parametrized solution snapshots stored in the POD reduced form, the interpolation process by the Bi-CITSGM is carried out in three main steps. First, the untrained POD eigenvalues are approximated by using the spline cubic interpolation. Then, the spatial and temporal POD bases are predicted by using the ITSGM method. Finally, two orthogonal calibration matrices are introduced in order to ensure the best match between the interpolated bases and their corresponding eigenvalues. These matrices are found as analytical solutions of two constrained optimization problems.   
In this paper, we propose a hyper reduced version referred to as Hyper Bi-CITSGM, that aim to accelerate the usual Bi-CITSGM by bringing the interpolation process to a significantly lower dimension. Based on this result, we design a reduced genetic algorithm (RGA) in which the full Navier-Stokes solver is replaced by the Hyper Bi-CITSGM method and where chromosomes are further enriched with additional virtual control parameters serving to ensure an optimal performance of RGA. The potential of RGA in terms of accuracy and CPU time is demonstrated on control problems of the flow past a cylinder and flow in a lid driven cavity when the Reynolds number value varies. 
\\
The paper is organized in the following manner. In the next section the Bi-CITSGM method \cite{OulghelouBiCITSGM2019Arxiv} proposed for nonlinear snapshots interpolation is briefly reviewed. In Section 3, the Hyper reduced approach of the Bi-CITSGM is introduced. Elements of the GA used for this study (RGA) and results on the control problems of the flow past a cylinder and a lid driven cavity are provided respectively in Section 4 and Section 5. Finally, the paper is concluded with a summary in Section 6.

%% file: Tex_Files/intro_Review_Bi-ITSGM.tex
Let $\trainingparam{i}\in\mathbb{R}^{\dimParam}$, $i = 1 ,\dots, \nbrParam$, be a set of parameters and $\trainedSnapMat{i}$ the associated parametrized snapshots matrices whose columns are solutions at different time instants of a non linear physical problem.
Each $\trainedSnapMat{i}$ is a $\dimVecorsSnap\times \nbrSnap$ matrix, where $\dimVecorsSnap$ is the number of spatial degrees of freedom and $\nbrSnap$ the number of time instants.
A classically frequented question is : using the set of matrices $\trainedSnapMat{i}$, is it possible to efficiently predict $\InterpSnapMat$ the matrix of snapshots for a new untrained parameter $\untrainedparam \neq \trainingparam{i}$. Unfortunately, this task is not straight forward given the nonlinear dependency of solutions to parameters. In such cases, usual interpolation techniques fail in general. To overcome this issue, it is possible to use the sophisticated interpolation approach Bi-CITSGM \cite{OulghelouBiCITSGM2019Arxiv} designed specifically for non linear snapshots matrices interpolation. This method is based on ITSGM method introduced for POD bases interpolation \cite{Amsallem,OULGHELOUAMC2018,Oulghelou2017}.
%
%
In the next subsections, a brief overview of the ITSGM and Bi-CITSGM methods is given. Further details can be found in \cite{ OulghelouBiCITSGM2019Arxiv,Amsallem}.

%% file: Tex_Files/Subspaces_interpolation.tex
\subsection{Geodesic Exponential and Logarithmic maps}
The ITSGM method proposed by Amsallem and Farhat \cite{Amsallem} 
is based on differential geometry tools involving geodesic Exponential and Logarithmic mappings in the Grassmann manifold. The Grassmann manifold $\GrassmanifoldSpace{}$ is defined as the set of all $\nbrModes$-dimensional subspaces in $\mbb{R}^{N}$, $0 \leq \nbrModes \leq N$. A point  $[\Xbasis{}] \in \GrassmanifoldSpace{}$ can be defined by the equivalence class \cite{Edelman1998, Boumal-2015}   
$$ [\Xbasis{}] = \{ \Xbasis{} \MatCalib{Q} \mid  \ \  \MatCalib{Q}\in\OrthGroupe{\nbrModes}\}$$ 
where $\Xbasis{}$ is a $N$ by $\nbrModes$ matrix with orthogonal columns, i.e, $\Xbasis{}^T\Xbasis{} = \Mat{I}_{\nbrModes}$, and $\OrthGroupe{\nbrModes}$ is the group of all $\nbrModes\times\nbrModes$ orthogonal matrices.
The geodesic distance $\distGrass{\Xbasis{}}{\XbasisEnd{}}$ between two points $[\Xbasis{}]$ and $[\XbasisEnd{}]$ in the Grassmann manifold is defined as the summation of squared principal angles \cite{Boumal-2015}
\begin{equation}
\distGrass{\Xbasis{}}{\XbasisEnd{}} = \sqrt{\displaystyle{\overset{}{\underset{i}{\sum}}} \arccos^2(\sigma_i)} 
\end{equation} 
where $\sigma_i$ are the singular values of $\Xbasis{}^T \XbasisEnd{}$. 
At each point $[\Xbasis{}]$ of the manifold $\GrassmanifoldSpace{}$, there exists a tangent space $\TangsubSpaceGrass{\Xbasis{}}{}$ of the same dimension \cite{Edelman1998, Absil} and a unique geodesic path \footnote{A geodesic between two points of the Grassmann manifold is the  path that minimizes the geodesic distance \cite{Wald}.} starting from $[\Xbasis{}]$ in every direction $\InitVel{} \in \TangsubSpaceGrass{\Xbasis{}}{}$, giving us the exponential map $\Exp{\Xbasis{}} : \TangsubSpaceGrass{\Xbasis{}}{} \longrightarrow \GrassmanifoldSpace{}$. Let $U \Sigma V^T$ be the thin SVD of the initial velocity $\InitVel{}$, the exponential of $\InitVel{}$ is given by
\begin{equation}\label{Eq:Exp_map}
[\XbasisEnd{}] = \txt{span} \{ \Xbasis{} V \cos(\Sigma) + U \sin(\Sigma)  \}
\end{equation}
Let us denote $\Log{\Xbasis{}}$ the inverse map to $\Exp{\Xbasis{}}$, which is defined only in a certain neighbourhood of $[\Xbasis{}]$. If $\Exp{\Xbasis{}}(\InitVel{}) = [\XbasisEnd{}]$, then $\InitVel{}$ is the vector determined as follows
\begin{equation}\label{Eq:Log_map}
\InitVel{} = \Log{\Xbasis{}}([\XbasisEnd{}]) = U \arctan(\Sigma) V^T
\end{equation}
where $U\Sigma V^T$ is the thin SVD of $(I-\Xbasis{} \Xbasis{}^T)\XbasisEnd{}(\Xbasis{}^T\XbasisEnd{})^{- 1}$ and $\Log{\Xbasis{}}([\Xbasis{}]) = 0$.
\subsection{ITSGM method}
Let $\Sbase{1}, \Sbase{2}, \dots, \Sbase{\nbrParam}$ be a set of parametrized POD \footnote{The POD appoach is not described in this article. For more details, the reader is referred to \cite{Sirovich}.} bases and $[\Sbase{1}], [\Sbase{2}], \dots, [\Sbase{\nbrParam}]$ the associated subspaces belonging to the Grassmann manifold.
The ITSGM problem is announced as follows : by using the definition of geodesic paths, Exponential and Logarithmic mappings, find an approximation of the subspace $[\InterpSbase]$ corresponding to a new untrained parameter $\untrainedparam \neq \ctrlvar_i$.
The first step of the ITSGM method is to chose arbitrarily a reference point $[\Sbase{i_0}]$ where $ i_0\in\{1,\dots,\nbrParam\}$. This reference point \footnote{New recently developed methods that do not require a reference can be found in \cite{TheseRolando, Rolando_IDW_article}.} is considered as the starting of the geodesic paths linking $[\Sbase{i_0}]$ to the rest of sampling subspaces $[\Sbase{i}]$. Now, by using the Logarithmic mapping, the initial velocity $\InitVel{i}$ for each geodesic path starting from $[\Sbase{i_0}]$ and ending at $[\Sbase{i}]$ can be calculated. Given that the Tangent space $\TangsubSpaceGrass{\Xbasis{i_0}}{}$ is a flat space, standard interpolation techniques such as Lagrange, RBF, Spline ...etc can be used. As a result,the initial velocity $\InterpGrassVel$ of the geodesic path linking the point $[\Sbase{i_0}]$ to the point $[\InterpSbase]$ can be interpolated. Finally, by using the geodesic Exponential mapping, an approximation of the subspace $[\InterpSbase]$ can be found. The steps of ITSGM are summarized in algorithm \ref{Alg:ITSGM}.
\\
\begin{algorithm}[H]
\begin{itemize}
\item[\hspace{10pt}]
\begin{itemize}
\item[ \textit{step 1}] Choose the origin point of tangency, for example $[\Sbase{i_0}]$ where $i_0 \in \{1,\dots,\nbrParam\}$.
\item[ \textit{step 2}] For $i \in \{1,\dots,\nbrParam\}$, map the point $[\Sbase{i}] \in \GrassmanifoldSpace{}$ to $\InitVel{i} \in \TangsubSpaceGrass{\Xbasis{i_0}}{}$ such that $\InitVel{i} = \Log{\Sbase{i_0}}(\Sbase{i})$ is the vector represented by 
$$ \InitVel{i} = U_i \arctan(\Sigma_i) V_i^T $$
where $U_i \Sigma_i V_i^T = (I-\Sbase{i_0} \Sbase{i_0}^T)\Sbase{i}(\Sbase{i_0}^T\Sbase{i})^{- 1}, i=1, \dots, \nbrParam$, are thin SVD.
\item[ \textit{step 3}] Interpolate the initial velocities $\InitVel{1}, \InitVel{2}, \dots, \InitVel{\nbrParam}$ for the untrained parameter $\untrainedparam$ using a standard interpolation and obtain $\InterpGrassVel$.
\item[ \textit{step 4}] Finally by the exponential mapping, map the interpolated velocity $\InterpGrassVel$ back to the Grassmann manifold. The matrix representation of the interpolated subspace is given by
$$
    \InterpSbase = \Sbase{i_0} \widetilde{V} \cos(\widetilde{\Sigma}) + \widetilde{U} \sin(\widetilde{\Sigma})
$$
where $\widetilde{U} \widetilde{\Sigma} \widetilde{V}^T$ is the thin SVD of the initial velocity vector $\InterpGrassVel$.
\end{itemize}
\end{itemize}
\caption{ITSGM}
\label{Alg:ITSGM}
\end{algorithm}

%% file: Tex_Files/Enhanced_Bi-CITSGM.tex
Let $\trainedSnapMat{i}\in\mbb{R}^{\dimVecorsSnap\times \nbrSnap}$, $i = 1 ,\dots, \nbrParam$, be a set of parametrized snapshots matrices \footnote{The parametrized snapshots matrices can result from CFD calculations or from experimental data.} whose columns are the solutions of a non linear physical problem.
The aim of the Bi-CITSGM is to appropriately use the existing matrices $\trainedSnapMat{i}$ in order to approximate $\InterpSnapMat$ for a new parameter $\untrainedparam \neq \trainingparam{i}$. 
%
The offline stage of the Bi-CITSGM method assumes that for each parameter $\trainingparam{i}$, the matrix $\trainedSnapMat{i}$ is approximated in a POD basis\footnote{The POD basis is calculated by using the Euclidean inner product.} of dimension $\nbrModes$  as follows
\begin{equation}\label{POD approx sampling}
\trainedSnapMat{i} \approx \Sbase{i} \trainedSingVal{i} \Tbase{i}^T, \hspace*{0.3cm} i=1,\dots,\nbrParam,
\end{equation}
where $\Sbase{i}\in \mbb{R}^{\dimVecorsSnap\times \nbrModes}$ and $\Tbase{i} \in \mbb{R}^{\nbrSnap\times \nbrModes}$ are respectively the left and right singular vectors of $\trainedSnapMat{i}$, and $\trainedSingVal{i}\in \mbb{R}^{\nbrModes\times \nbrModes}$ the corresponding matrix of singular values.
%
%
%
In the online stage, the matrix of singular values $\InterpSingVal$ is first approximated by using spline cubic interpolation.
Next, the spatial and temporal bases $\InterpSbase$ and $\InterpTbase$ are respectively approximated by interpolating $\Sbase{i}$ and $\Tbase{i}$ using the ITSGM method. In order to ensure a well orientation of modes in the POD sampling bases, the signs are adjusted such that 
%
the $j^{\txt{th}}$ spatial and temporal modes $\Sbase{k}^j$ and $\Tbase{k}^j$ are multiplied by $-1$ if the following condition is fulfilled
$$ || \Sbase{k_0}^j - \Sbase{k}^j ||_2 > || \Sbase{k_0}^j + \Sbase{k}^j ||_2$$
where $k_0$ is the index of the reference basis determined as
$$ k_0 = \underset{i\in\{1,\dots,\nbrParam\}}{\txt{\textbf{argmin}}} \ \ \distGrass{\InterpSbase}{\Sbase{i}} $$
%
%
Finally the interpolated bases are calibrated by two orthogonal matrices obtained as analytical solutions of the following two optimization problems
\begin{equation}\label{Eq : space_minimization_problem_Bi-CITSGM}
\underset{\CalibMatrix{Q}{x}\in \mcal{O}(q)}{\txt{\textbf{argmin}}} \ \ \somme{i}{1}{\nbrParam} \omega_i  || \InterpSbase  \CalibMatrix{Q}{x} - \Sbase{i}||_F^2
\hspace*{1cm}
\underset{\CalibMatrix{Q}{t}\in \mcal{O}(q)}{\txt{\textbf{argmin}}} \ \ \somme{i}{1}{\nbrParam} \kappa_i  || \InterpTbase  \CalibMatrix{Q}{t} - \Tbase{i}||_F^2
\end{equation}
where $||\cdot||_F$ is the Frobenius norm and $\omega_i$ and $\kappa_i$ are the GIDW (Grassmann Inverse Distance weighting) weights given for $m,l>1$ by 
\begin{equation}\label{Bi-CITSGM_weights}
\omega_i = \myfrac{\distGrass{\InterpSbase}{\Sbase{i}}^{-m}}{\somme{k}{1}{\nbrParam} \distGrass{\InterpSbase}{\Sbase{k}}^{-m}} 
\hspace*{1.5cm}
\kappa_i = \myfrac{\distGrass{\InterpTbase}{\Tbase{i}}^{-l}}{\somme{k}{1}{\nbrParam} \distGrass{\InterpTbase}{\Tbase{k}}^{-l}} 
\end{equation}
Let $\Mat{M}_x$ and $\Mat{M}_t$ be the following $\nbrModes\times\nbrModes$ matrices 
$$\Mat{M}_x = \InterpSbase^T \somme{i}{1}{\nbrParam} \omega_i \Sbase{i} 
\hspace*{1cm}
\Mat{M}_t = \InterpTbase^T \somme{i}{1}{\nbrParam} \kappa_i \Tbase{i}$$
The analytical solutions of optimization problems \eqref{Eq : space_minimization_problem_Bi-CITSGM} writes as follows
$$\CalibMatrix{Q}{x} =  \Mat{\xi} \Mat{I}_{_{\CalibMatrix{Q}{x}}}^+ \Mat{\eta}^T
\hspace*{1cm}
\CalibMatrix{Q}{t} =  \Mat{\zeta} \Mat{I}_{_{\CalibMatrix{Q}{t}}}^+ \Mat{\rho}^T$$
where $\Mat{\xi}$ and $\Mat{\eta}$ (resp.  $\Mat{\zeta}$ and $\Mat{\rho}$) are the left and right singular vectors of $\Mat{M}_x$ (resp. $\Mat{M}_t$) and $\Mat{I}_{_{\CalibMatrix{Q}{x}}}^+$ (resp. $\Mat{I}_{_{\CalibMatrix{Q}{t}}}^+$) is the diagonal matrix whose elements are equal to $1$ for non zero singular values and $0$ otherwise.
The steps of the Bi-CITSGM are summarized in algorithm \ref{Alg:Bi-CITSGM}. For a detailed review of the method, the interested reader is referred to \cite{OulghelouBiCITSGM2019Arxiv}. 
\begin{algorithm}[H]
\begin{itemize}
\item[]
\begin{itemize}
\item[\textit{\textbf{Offline :}}]
\item[\textit{step 1 :}] Perform POD of order $\nbrModes$ of the sampling snapshots matrices
$$ \trainedSnapMat{i} \approx \Sbase{i} \trainedSingVal{i} \Tbase{i}^T, \hspace*{0.3cm} i=1,\dots,\nbrParam $$
\item[\textit{\textbf{Online :}}]
\item[\textit{step 2 :}] Using spline cubic, interpolate $\trainedSingVal{i}$, $i=1,\dots,\nbrParam$, to obtain $\InterpSingVal$
\item[\textit{step 3 :}] Interpolate $[\Sbase{i}]$ and $[\Tbase{i}]$, $i=1,\dots,\nbrParam$, by using the ITSGM method (algorithm \ref{Alg:ITSGM}) and obtain the spatial and temporal bases $\InterpSbase$ and $\InterpTbase$
\item[\textit{step 4 :}] Find $ k_0 = \underset{i\in\{1,\dots,\nbrParam\}}{\txt{\textbf{argmin}}} \ \ \distGrass{\InterpSbase}{\Sbase{i}} $ and adjust sampling bases modes signs 
\item[\textit{step 5 :}] Calculate the weights $\omega_i$ and $\kappa_i$ using equations \eqref{Bi-CITSGM_weights}
\item[\textit{step 6 :}] Perform SVD decompositions 
$$\InterpSbase^T \somme{i}{1}{\nbrParam} \omega_i \Sbase{i} = \Mat{\xi} \Mat{\Theta} \Mat{\eta}^T
\hspace*{1cm}
\InterpTbase^T \somme{i}{1}{\nbrParam} \kappa_i \Tbase{i} = \Mat{\zeta} \Mat{\delta} \Mat{\rho}^T$$
\item[\textit{step 7 :}] Evaluate the calibration matrices $\CalibMatrix{Q}{x}$ and $\CalibMatrix{Q}{t}$ as follows
$$\CalibMatrix{Q}{x} =  \Mat{\xi} \Mat{I}_{_{\CalibMatrix{Q}{x}}}^+ \Mat{\eta}^T
\hspace*{1cm}
\CalibMatrix{Q}{t} =  \Mat{\zeta} \Mat{I}_{_{\CalibMatrix{Q}{t}}}^+ \Mat{\rho}^T$$
where $\Mat{I}_{_{\CalibMatrix{Q}{x}}}^+$ and $\Mat{I}_{_{\CalibMatrix{Q}{t}}}^+$ are diagonal matrices whose elements are equal to $1$ for non zero singular values and $0$ otherwise
\item[\textit{step 8 :}] Reconstruction of the interpolated snapshots matrix 
$$\InterpSnapMat = \InterpSbase \CalibMatrix{Q}{x} \InterpSingVal \CalibMatrix{Q}{t}^T \InterpTbase^T$$ 
\end{itemize}
\end{itemize}
\caption{Bi-CITSGM}
\label{Alg:Bi-CITSGM}
\end{algorithm}

%% file: Tex_Files/Hyper_Bi-CITSGM.tex
Recall that $\trainedSnapMat{i}\in\mbb{R}^{\dimVecorsSnap\times \nbrSnap}$, $i = 1 ,\dots, \nbrParam$, is the set of parametrized snapshots matrices whose columns are the solutions of a non linear physical problem. Each matrix $\trainedSnapMat{i}$ is approximated in a POD basis of dimension $\nbrModes$ as follows
\begin{equation}\label{POD approx sampling}
\trainedSnapMat{i} \approx \Sbase{i} \trainedSingVal{i} \Tbase{i}^T
\end{equation}
where $\Sbase{i}\in \mbb{R}^{\dimVecorsSnap\times \nbrModes}$ and $\Tbase{i} \in \mbb{R}^{\nbrSnap\times \nbrModes}$ are respectively the left and right singular vectors of $\trainedSnapMat{i}$, and $\trainedSingVal{i}\in \mbb{R}^{\nbrModes\times \nbrModes}$ the corresponding matrix of singular values.
The aim of the following section is to drastically reduce the dimensionality of the Bi-CITSGM problem. 
Consider the POD respectively of order $r$ and $s$, $r,s \leq \nbrModes \nbrParam$, of the following column block matrices 
\begin{equation*}
\begin{bmatrix} \Sbase{1} &  \Sbase{2} & \cdots & \Sbase{\nbrBases} \end{bmatrix} = \bm{\Phi}  \varrho \Mat{W}^T
\ \ \ \ \ \txt{and} \ \ \ \ \ 
\begin{bmatrix} \Tbase{1} & \Tbase{2} & \cdots & \Tbase{\nbrBases} \end{bmatrix} = \bm{\Lambda}  \Theta \Mat{Z}^T
\end{equation*}
where $\bm{\Phi}\in \mbb{R}^{\dimVecorsSnap \times r}$, $\Mat{W}\in \mbb{R}^{\nbrModes \nbrBases \times r}$, $\bm{\Lambda}\in \mbb{R}^{\nbrSnap \times s}$, $\Mat{Z}\in \mbb{R}^{\nbrModes \nbrBases \times s}$, $\varrho = diag\left(\rho_{1}, \rho_{2}, \cdots,\rho_{r}\right)$ and $\Theta = diag\left(\theta_{1}, \theta_{2}, \cdots,\theta_{s}\right)$.
Let $W_i\in \mbb{R}^{\nbrModes \times r}$ and $Z_i\in \mbb{R}^{\nbrModes \times s}$, $i = 1 ,\dots, \nbrParam$, be the row blocks of the matrices $\Mat{W}$ and $\Mat{Z}$, i.e.,
\begin{equation*}
\Mat{W} = \begin{bmatrix} W_1 \\ W_2 \\ \vdots \\ W_{\nbrBases} \end{bmatrix} \ \ \ \ \txt{and} \ \ \ \ 
\Mat{Z} = \begin{bmatrix} Z_1 \\ Z_2 \\ \vdots \\ Z_{\nbrBases} \end{bmatrix}
\end{equation*}
The trained snapshot matrices $\trainedSnapMat{i}$ can be expressed in the following manner
\begin{equation}\label{HNIMR snapshots matrix}
\trainedSnapMat{i} \approx \bm{\Phi}  \varrho  \trainedredSnapMat{i} \Theta \bm{\Lambda}^T
\end{equation}
where
$$\trainedredSnapMat{i} = W_i^T \Mat{\Sigma}_i Z_i$$
From expression \ref{HNIMR snapshots matrix}, it's obvious that $\bm{\Phi}$, $\varrho$, $\Theta$ and $\bm{\Lambda}$ remain constant, and that the matrix $\trainedredSnapMat{i}$ depends on the parameter $\ctrlvar_i$.
Thus, a new ensemble of reduced parametrized snapshot matrices can be generated. The main advantage of this new ensemble is that its elements $\trainedredSnapMat{i}$ are matrices of reduced size $r\times s$, where $r,s\leq \nbrModes \nbrParam\ll \dimVecorsSnap$.
As a result, instead of applying the Bi-CITSGM directly to the set of $\trainedSnapMat{i}$ involving manipulation of large matrices, it is more convenient in terms of memory storage and computational time to use the decomposition \eqref{HNIMR snapshots matrix} and apply the Bi-CITSGM to the new set of reduced snapshots matrices $\trainedredSnapMat{i}$. This hyper reduced version is referred to as Hyper Bi-CITSGM.
%
It's worth noting that if the order of truncation $r$ and $s$ are such as $r = s = \nbrParam \nbrModes$, the outputs of the Hyper Bi-CITSGM and the usual Bi-CITSGM are exactly identical. In the case $r,s < \nbrParam \nbrModes$, the outputs can be close to each others provided that $r$ and $s$ are properly chosen. In practice, $r$ and $s$ are chosen such that $\theta_r<\epsilon_r$ and $\rho_s < \epsilon_s$, where $\epsilon_r$ and $\epsilon_s$ are sufficiently small thresholds chosen by the user.
The steps of the Hyper Bi-CITSGM are listed in algorithm \ref{Alg:Hyper Bi-CITSGM}.
%
%
\\
\begin{algorithm}[H]
\begin{itemize}
\item[]
\begin{itemize}
\item[\textit{\textbf{Offline :}}]
\item[\textit{step 1 :}] Perform POD of order $\nbrModes$ of the sampling snapshots matrices
$$ \trainedSnapMat{i} \approx \Sbase{i} \trainedSingVal{i} \Tbase{i}^T, \hspace*{0.3cm} i=1,\dots,\nbrParam $$
\item[\textit{step 2 :}] Perform POD of order $r$ and $s$ ($r,s < \nbrParam \nbrModes$) respectively of the column block matrices 
\begin{align*}
\begin{bmatrix} \Sbase{1} &  \Sbase{2} & \cdots & \Sbase{\nbrBases} \end{bmatrix} &= \bm{\Phi}  \varrho \Mat{W}^T
\\
\begin{bmatrix} \Tbase{1} & \Tbase{2} & \cdots & \Tbase{\nbrBases} \end{bmatrix} &= \bm{\Lambda}  \Theta \Mat{Z}^T
\end{align*}
\item[\textit{step 3 :}] Extract the $\nbrModes \times r$ (resp. $\nbrModes \times s$) row block matrices $W_i$ (resp. $Z_i$) from $\Mat{W}$ (resp. $\Mat{Z}$), i.e.,   
$$
\Mat{W} = \begin{bmatrix} W_1 \\ W_2 \\ \vdots \\ W_{\nbrBases} \end{bmatrix} 
\ \ \ \ \ \txt{and} \ \ \ \ \ 
\Mat{Z} = \begin{bmatrix} Z_1 \\ Z_2 \\ \vdots \\ Z_{\nbrBases} \end{bmatrix}
$$
%
\item[\textit{step 3 :}] Construct the $r\times s$ reduced snapshots matrices $\trainedredSnapMat{i}$ as follows
$$ \trainedredSnapMat{i} = W_i^T \Mat{\Sigma}_i Z_i$$
\item[\textit{step 4 :}] Perform POD of order $\nbrModes$ of the sampling reduced snapshots matrices
$$ \trainedredSnapMat{i} \approx \phi_{i} \chi_{i} \alpha_{i}^T, \hspace*{0.3cm} i=1,\dots,\nbrParam $$
\item[\textit{\textbf{Online :}}]
\item[\textit{step 5 :}] Apply \textit{step 2} to \textit{step 8} of algorithm \ref{Alg:Bi-CITSGM} to the ensemble of matrices $\phi_{i}$, $\chi_{i}$ and $\alpha_{i}$ and obtain an approximation $\InterpRedSnapMat$ of the snapshot matrix associated to the new untrained parameter $\untrainedparam$
\item[\textit{step 6 :}] Reconstruct the interpolated snapshots matrix 
$$\InterpSnapMat = \bm{\Phi}\varrho  \InterpRedSnapMat  \Theta\bm{\Lambda}^T$$ 
\end{itemize}
\end{itemize}
\caption{Hyper Bi-CITSGM}
\label{Alg:Hyper Bi-CITSGM}
\end{algorithm}

%% file: Tex_Files/Comp_complexity_Hyper_Bi-CITSGM.tex
It was shown that in the case of univariate interpolation, the computational complexity of the Bi-CITSGM is proportional to $\mcal{O}(\dimVecorsSnap \nbrModes^2)$ \cite{Amsallem, OulghelouBiCITSGM2019Arxiv}, where $\dimVecorsSnap$ is the number of degrees of freedom of the parametrized physical problem and $\nbrModes$ is the dimension of the spatial POD basis used to approximate its solutions. 
%
By using the Hyper Bi-CITSGM method, the computational complexity reduce to be proportional to $\mcal{O}(\nbrParam \nbrModes^3)$, where $\nbrParam$ is the number of trained parameters. 
This suggests that the proposed approach is computationally efficient.

%% file: Tex_Files/ctrl_prob.tex
Consider the constrained nonlinear optimization problem of the form
\begin{equation}\label{ctrl_NS}
\underset{\ctrlvar}{\min} \ \ \functional( \bm{y}(\ctrlvar) , \ctrlvar)
\hspace*{0.5cm}
\txt{subject to }
\hspace*{0.5cm}
\constraint(\bm{y}(\ctrlvar), \ctrlvar) = 0
\end{equation}
where $\functional$ is the cost function and $\bm{y}$ and $\ctrlvar$ denote the state and control variables related to each other trough the constraint mapping $\constraint$. In flow optimal control, $\bm{y}$ may correspond to velocity, pressure or temperature of the fluid, $\ctrlvar$ to the Reynolds number, Strouhal number or the angle of attack, and the mapping $\constraint$ to the Navier-Stokes equations or the experimental protocol. In the following, we are interested in solving flow control problems of type \eqref{ctrl_NS} using genetic algorithms. First, the principle of genetic algorithms is introduced, then a reduced genetic algorithm (RGA) approach for near real time optimal control is proposed. 

%% file: Tex_Files/GA_overview.tex
\begin{figure}[hbtp!]
\centering
\begin{tikzpicture}
\node[rectangle,draw=black, ultra thick] (a) {\textbf{initialization}};
\node[rectangle,draw=black, ultra thick, below=0.5cm of a] (b) {\textbf{stop criterion}};
\node[rectangle,draw=black, ultra thick, right=1cm of b] (g) {\textbf{termination}};
\node[rectangle,draw=black, ultra thick, below=0.5cm of b] (c) {\textbf{selection}};
\node[rectangle,draw=black, ultra thick, below=0.3cm of c] (d) {\textbf{crossover}};
\node[rectangle,draw=black, ultra thick, below=0.3cm of d] (e) {\textbf{mutation}};

\node[right =0.35cm of b] (yp) {};
\node[thick, above =-0.2cm of yp] (y) {yes};

\node[below =0.1cm of b] (np) {};
\node[thick, right =-0.15cm of np] (n) {no!};

\draw[ultra thick, ->] (a) to (b) ;
\draw[ultra thick, ->] (b) to (c) ;
\draw[ultra thick, ->] (c) to (d) ;
\draw[ultra thick, ->] (d) to (e) ;
\draw[ultra thick, ->] (b) to (g) ;
\draw [-, rounded corners, ultra thick] (0,-0.5) -- (-2.25,-0.5) -- (-2.25,-4.5) -- (0,-4.5)--(0,-4);
\end{tikzpicture}
\caption{Outline of the canonical GA used for optimization problems.}
\label{fig:GA}
\end{figure}
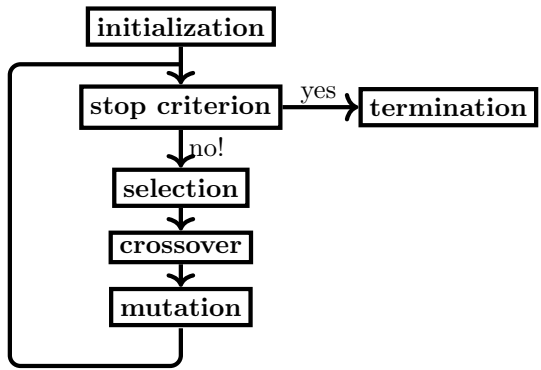

%
%
%
%
%
%
%
GA starts with a random set (of size $N_{_{\txt{chrom}}}$) of chromosomes $\ctrlvar_1, \ctrlvar_2, \dots, \ctrlvar_{N_{_{\txt{chrom}}}} $. Each chromosome $\ctrlvar_j$ is evaluated using the objective function $\functional$ and constraint information (Navier-Stokes equations or experiment), and a fitness value $f(\ctrlvar_j)$ is assigned. 
Then, three main genetic operators (selection, crossover and mutation) modeled on the Darwinian concepts of natural selection and evolution are applied to the population in order to create a new hopefully better population.
In the selection, a new population of chromosomes is chosen to survive based on their fitness value. 
When the chromosome has larger fitness, it has higher probability of being reproduced and passed down into the next generation.
%
%
%
In the crossover, all surviving chromosomes are randomly paired. At a given crossover probability $P_c$, the pairs
exchange genes at a random locus. In practice, a random number ranging from $0$ to $1$ is generated. Then, if the random number is smaller than $P_c$, the crossover take place, and two new chromosomes are created to replace the original chromosomes. However, if the random number is greater than $P_c$, the two chromosomes in the original pair remain into the next generation.
In the mutation, genes in each chromosome are randomly altered at a mutation probability of $P_m$. Similar to the crossover, the
probability of the mutation is determined by a random number ranging from $0$ to $1$.
These genetic operations are repeated till a global optimal solution is approached. In our case, the GA is iterated until the number of renewed generations reached a predetermined number. The best recorded solution at the last generation is declared as the optimized solution. A sketch of a canonical genetic algorithm is shown in figure \ref{fig:GA}.
\\
Despite their robustness, a serious weakness of GAs in flow control is their high requirements in terms of CPU time and memory storage. In fact, we have to solve and store solutions of the Navier-Stokes equations many times for each population, which is numerically very expensive. To overcome this issue we use the Hyper Bi-CITSGM proposed in section 3 to design a reduced genetic algorithm strategy (RGA) capable of performing near real time flow control. The next subsection is dedicated to describe the proposed RGA.
%

%% file: Tex_Files/used_GA.tex
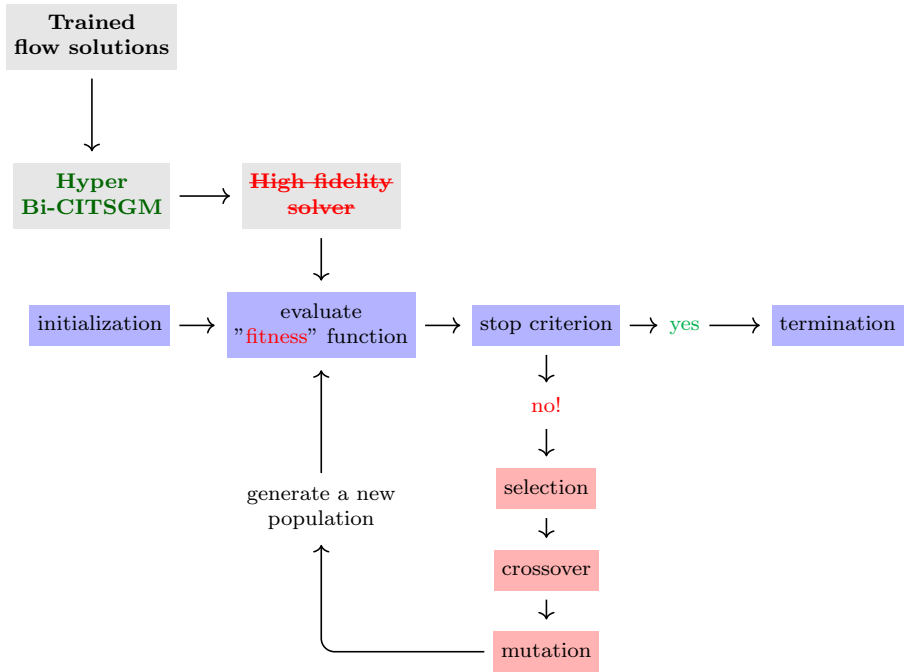
\begin{figure}[hbtp!]
\centering

\begin{tikzpicture}[
    pre/.style={<-,shorten <=1pt,semithick},
    post/.style={->,shorten >=1pt,semithick}
    ]

\node (n1) at (0.,0.) {\footnotesize \tikz[baseline]{\node[fill=blue!30,anchor=base] {\minibox[c]{
initialization  }};}};

\node[right = 0.5cm of n1] (n2)  {\footnotesize  \tikz[baseline]{\node[fill=blue!30,anchor=base] {\minibox[c]{
evaluate \\ "\textcolor{red}{fitness}" function }};}};

\node[above = 0.6cm of n2] (resNS) {\footnotesize \tikz[baseline]{\node[fill=gray!20,anchor=base] {\minibox[c]{
\textcolor{red}{\textbf{\sout{High fidelity}}} \\ \textcolor{red}{\textbf{\sout{solver}}}
 }};}};

\node[left = 0.7cm of resNS] (HNIMR) {\footnotesize \tikz[baseline]{\node[fill=gray!20,anchor=base] {\minibox[c]{
\textcolor{green!40!black}{\textbf{Hyper}} \\ \textcolor{green!40!black}{\textbf{Bi-CITSGM}} }};}};

\node[above = 1cm of HNIMR] (sampling) {\footnotesize \tikz[baseline]{\node[fill=gray!20,anchor=base] {\minibox[c]{
\textbf{Trained} \\ \textbf{flow solutions}  }};}};

\node[right = 0.5cm of n2] (n3)  {\footnotesize  \tikz[baseline]{\node[fill=blue!30,anchor=base] {\minibox[c]{
stop criterion  }};}};

\node[right = 0.4cm of n3] (oui) {\footnotesize  \minibox[c]{
\textcolor{green!70!blue}{yes}
}};

\node[below = 0.4cm of n3] (non) {\footnotesize  \minibox[c]{
\textcolor{red}{no!}
}};


\node[right = 0.7cm of oui] (fin)  {\footnotesize  \tikz[baseline]{\node[fill=blue!30,anchor=base] {\minibox[c]{
termination  }};}};

\node[below = 0.4cm of non] (g1) {\footnotesize  \tikz[baseline]{\node[fill=red!30,anchor=base] {\minibox[c]{
selection  }};}};

\node[below = 0.3cm of g1] (g2)  {\footnotesize  \tikz[baseline]{\node[fill=red!30,anchor=base] {\minibox[c]{
crossover  }};}};

\node[below = 0.3cm of g2] (g3)  {\footnotesize  \tikz[baseline]{\node[fill=red!30,anchor=base] {\minibox[c]{
mutation  }};}};

\node[below = 1.4cm of n2] (newpop) {\footnotesize  \minibox[c]{
generate a new 
\\
population
}};


\draw[post, ->] (HNIMR) to (resNS);
\draw[post, ->] (resNS) to (n2);

\draw[post, ->] (n1) to (n2);
\draw[post, ->] (n2) to (n3);
\draw[post] (n3) to (oui);
\draw[post, ->] (oui) to (fin);

\draw[post] (n3) to (non);
\draw[post, ->] (non) to (g1);
\draw[post, ->] (g1) to (g2);
\draw[post, ->] (g2) to (g3);

\draw[post,rounded corners=5pt] (g3.west) -| (newpop.south);
\draw[post] (newpop) to (n2);

\draw[post] (sampling) to (HNIMR);

\end{tikzpicture}
 \caption{Outline of the proposed RGA.}
\label{fig:RGA}
\end{figure}

In the present paper, a real coded Reduced Genetic Algorithm (RGA) is proposed. As shown in figure \ref{fig:RGA}, this algorithm consists of replacing the high fidelity solver in the stage of fitness evaluation by the cheap Hyper Bi-CITSGM interpolation approach. In order to enhance the performance of Hyper Bi-CITSGM inside the RGA, virtual genes are added to chromosomes. These are the GIDW (Grassmann Inverse distance weighting) powers $m,l>1$ appearing in the definition of wights $\omega_i$ and $\kappa_i$ in equations \eqref{Bi-CITSGM_weights}, and the time phase shift $\delta$. More particularly, the phase shift $\delta$ is introduced for tracking functional minimization problems. 
%
%
%
\\
Using a set of precalculated parametrized flow solutions for which steps \textit{1} to \textit{4} of algorithm \ref{Alg:Hyper Bi-CITSGM} were performed, the constrained nonlinear optimization problem solved by RGA is defined as follows
\begin{equation}\label{ctrl_NS}
\underset{\bar{\ctrlvar}}{\min} \ \ \functional( \bm{y}(\bar{\ctrlvar}), \bar{\ctrlvar} )
\hspace*{0.5cm}
\txt{subject to }
\hspace*{0.5cm}
\constraint_{\txt{red}}(\bm{y}(\bar{\ctrlvar}), \bar{\ctrlvar}) = 0
\end{equation}
where $\bm{y}(\bar{\ctrlvar})$ is the interpolated state solution at $\bar{\ctrlvar} = \{\ctrlvar, m, l, \delta \}$ through the constraint mapping $\constraint_{\txt{red}}$ described by steps \textit{5} and \textit{6} of algorithm \ref{Alg:Hyper Bi-CITSGM}.
In the RGA, the $j^{\txt{th}}$ chromosome corresponds to control candidate $\bar{\ctrlvar}_j = \{\ctrlvar_j, m_j, l_j, \delta_j \}$ where  $\ctrlvar_j$, $m_j$, $l_j$ and $\delta_j$ represent its genes.
The fitness function $f$ associated to the $j^{\txt{th}}$ chromosome is defined by the inverse of objective function
$$f(\bar{\ctrlvar}_j) = \myfrac{1}{\functional(\bm{y}(\bar{\ctrlvar}_j), \bar{\ctrlvar}_j)} $$
The probability of reproduction $P_s^j$ of the $j^{\txt{th}}$ chromosome is calculated as follows
\begin{equation}\label{selection_probability_formula}
P_s^j =  \myfrac{f(\bar{\ctrlvar}_j)}{\somme{i}{1}{N_{_{\txt{chrom}}}} f(\bar{\ctrlvar}_i)}
\end{equation}
Using this reproduction probability, $N_{_{\txt{chrom}}}$ solutions from the current generation are selected by the roulette rule \cite{Goldberg1989} to survive for the next generation. These reproduced solutions are afterwards modulated by the crossover and mutation operators.
%
%
%
A common difficulty arising in GA's implementation is how to determine the algorithm parameters, such as crossover and mutation probabilities. Generally, no specific criteria exists for such determinations. In this article, we conducted several trial runs in which we changed these parameters and chose a set of parameters that seemed to be appropriate with respect to convergence rate. For the proposed RGA, the mutation probability was set to $40\%$ and a simple one-point crossover operator was used with a probability of $60\%$. 
In the next subsections, the effectiveness of RGA is tested on the control problems of flow past a cylinder and flow in a lid driven cavity by acting on Reynolds number value.

%% file: Tex_Files/GA_Cylinder.tex
\begin{figure}[hbtp!]
\centering 
\includegraphics[width=\linewidth]{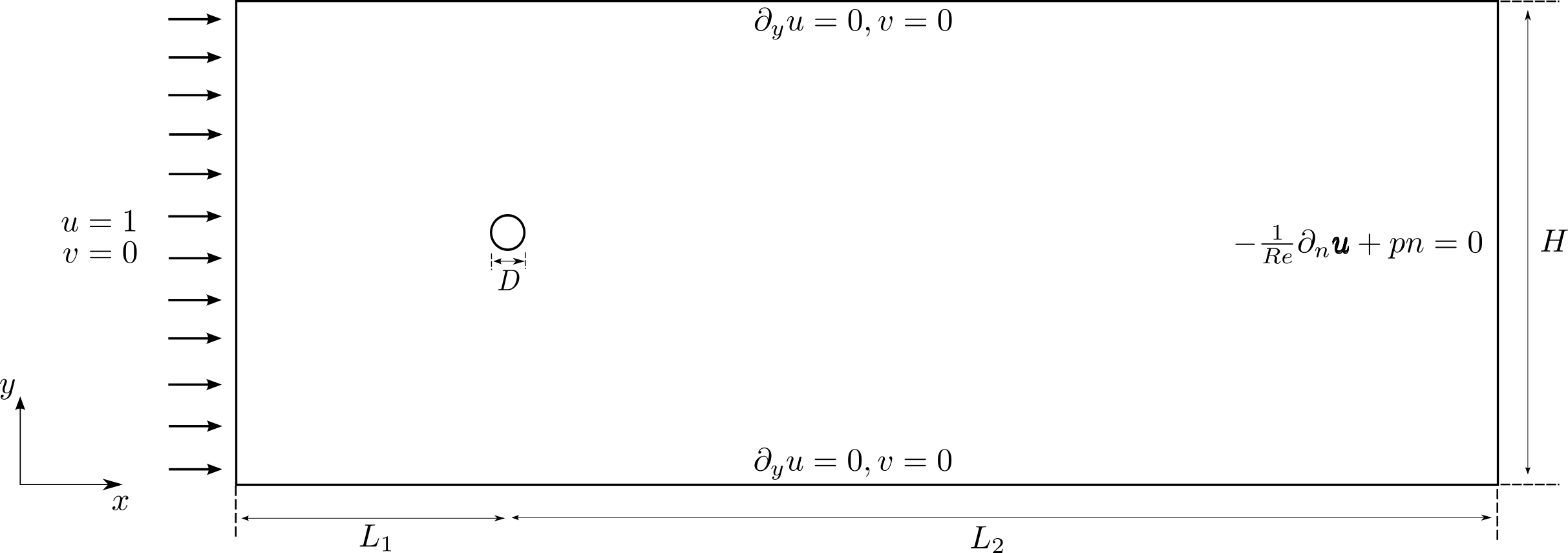}
\caption{Two-dimensional domain and boundary conditions for the problem of flow past a cylinder.}
\label{Fig:flow past cylinder}
\end{figure}
\begin{figure}[hbtp!]
\centering
\begin{subfigure}{0.5\textwidth}
\includegraphics[width=\linewidth]{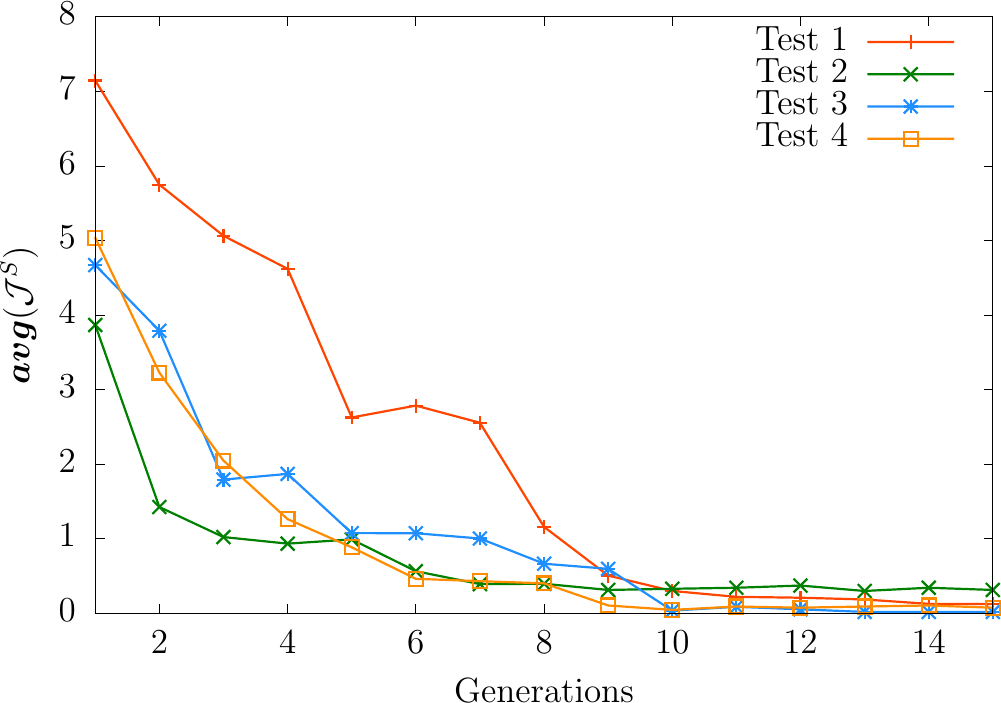}
\caption{minimization of $\functional^S$}
\end{subfigure}%
~ 
\begin{subfigure}{0.5\textwidth}
\includegraphics[width=\linewidth]{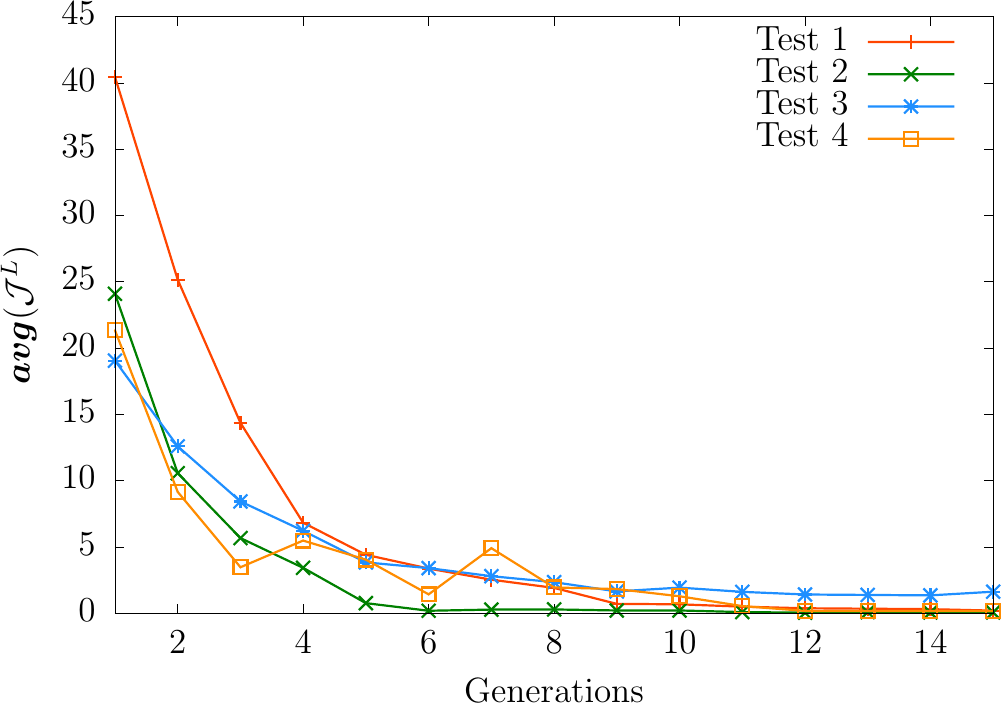}
\caption{minimization of $\functional^L$}
\end{subfigure}
\caption{Evolution of the averaged functionals $\bm{avg}(\functional^S)$ and $\bm{avg}(\functional^L)$ with respect to generations. }
\label{average_func_cylinder}
\end{figure}
\begin{figure}[hbtp!]
\begin{subfigure}{\textwidth}
\includegraphics[width=0.32\linewidth]{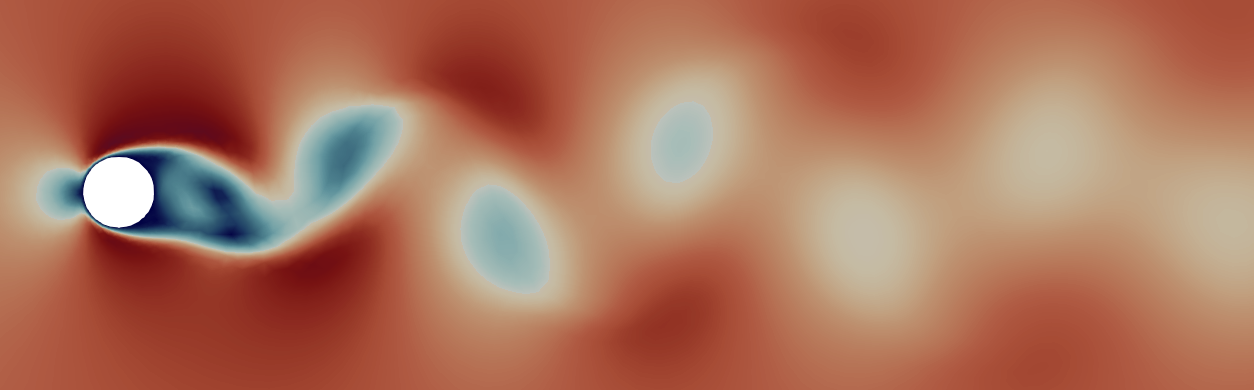}%
\hspace*{0.25cm}\includegraphics[width=0.32\linewidth]{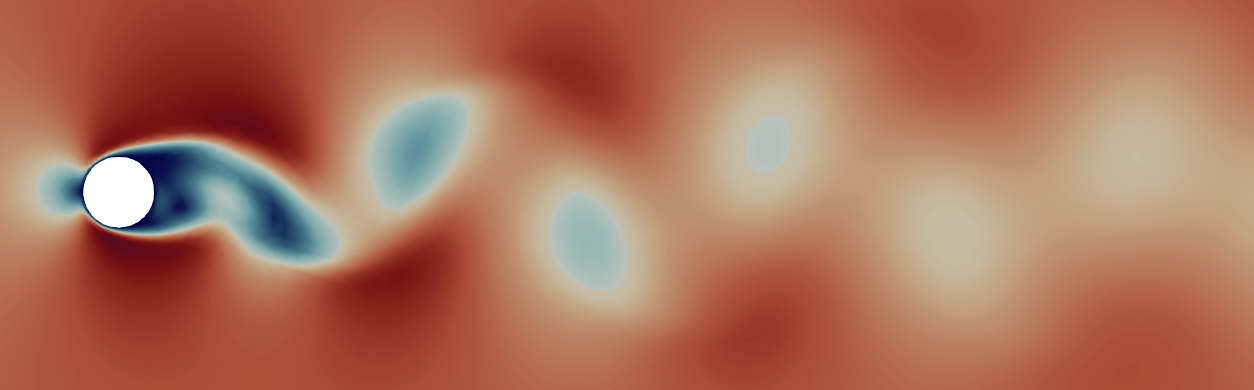}%
\hspace*{0.25cm}\includegraphics[width=0.32\linewidth]{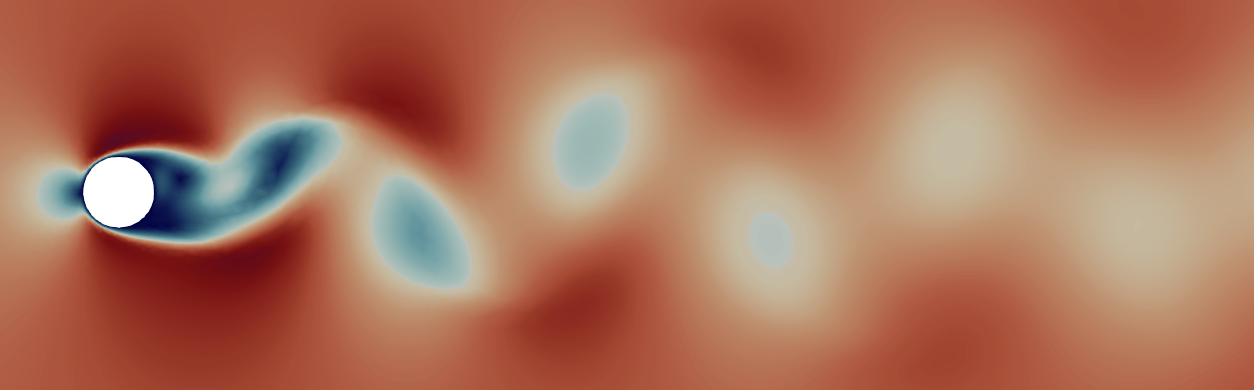}
\caption{Optimal flow solution }
\end{subfigure}
\begin{subfigure}{\textwidth}
\includegraphics[width=0.32\linewidth]{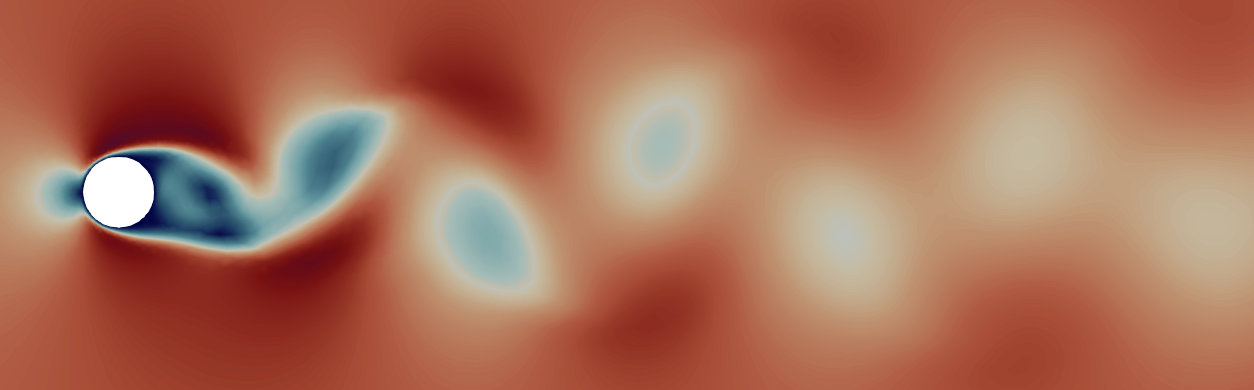}%
\hspace*{0.25cm}\includegraphics[width=0.32\linewidth]{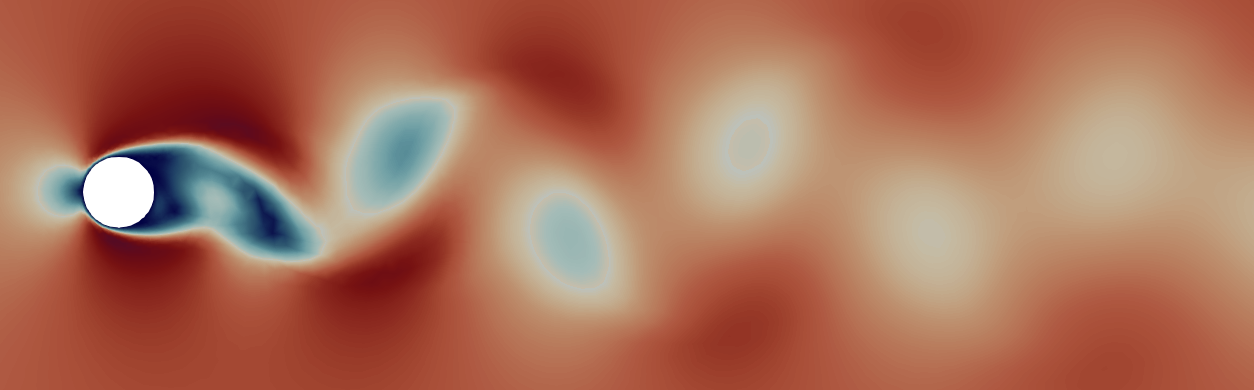}%
\hspace*{0.25cm}\includegraphics[width=0.32\linewidth]{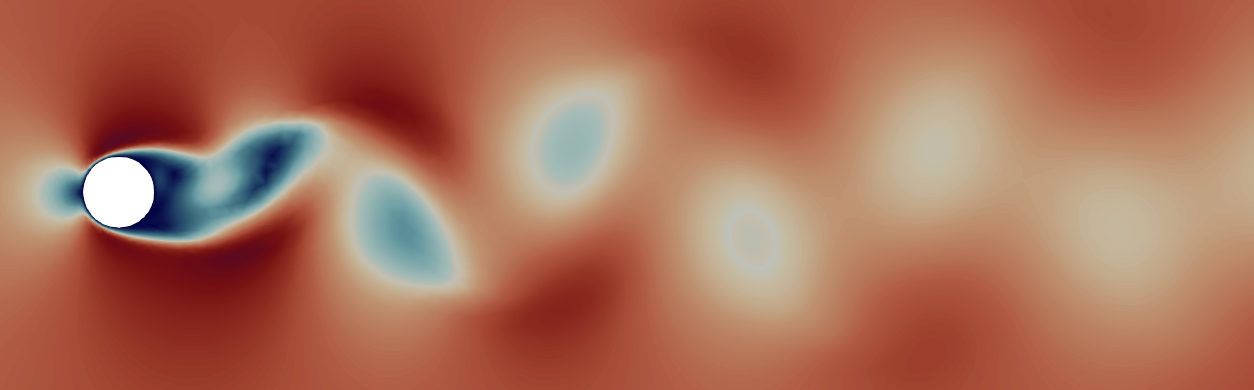}
\caption{flow solution obtained at generation $12$ of RGA applied to minimize $\functional^S$}
\end{subfigure}
\begin{subfigure}{\textwidth}
\includegraphics[width=0.32\linewidth]{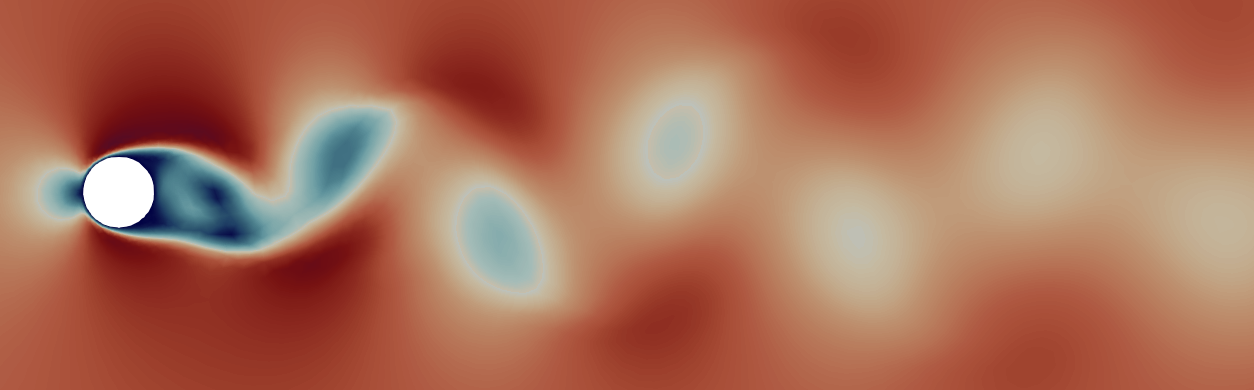}%
\hspace*{0.25cm}\includegraphics[width=0.32\linewidth]{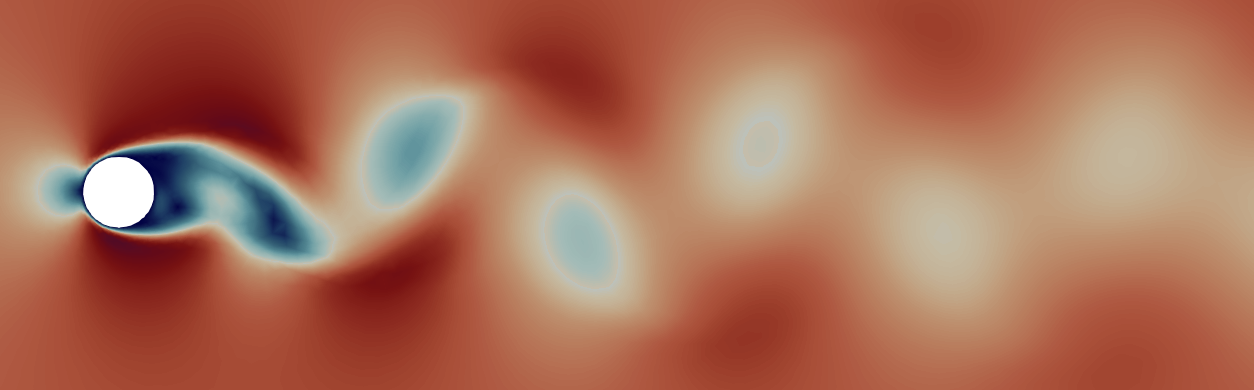}%
\hspace*{0.25cm}\includegraphics[width=0.32\linewidth]{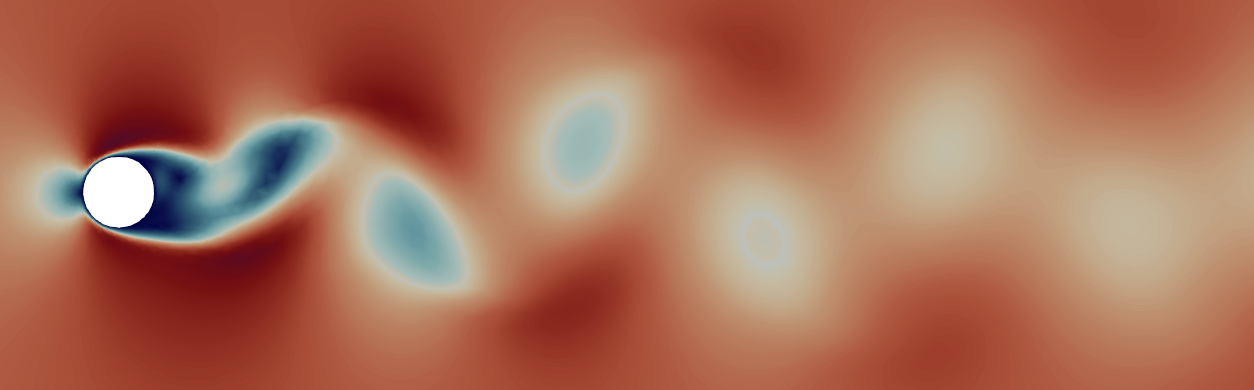}
\caption{flow solution obtained at generation $12$ of RGA applied to minimize $\functional^L$}
\end{subfigure}
\caption{Evolution of the target and predicted RGA solutions for Test 2 ($Re_{_{\txt{opt}}} = 160$) at tree instants $t=T/4$ (left), $t=T/2$ (middle) and $t=T$ (right). Here, $T$ represents a single flow period at $Re = 160$.}
\label{fig:3snapflow_RE160}
\end{figure}

Consider the two dimensional flow past a cylinder of diameter $D$ depicted in figure \ref{Fig:flow past cylinder}. 
The problem domain is rectangular with length $H = 30D$ and width  $45 D$ and contains a cylinder situated at $L_1 = 10D$ from the left boundary and $H/2$ from the lower boundary. The fluid dynamics of the flow is driven by an inlet velocity $U$ of a unit magnitude, which enters from the left boundary of the domain, and is allowed to flow past through the right boundary of the domain. Free slip boundary conditions are applied to the horizontal edges whilst no slip boundary condition are considered on the cylinder’s wall. The Reynolds number for this flow is given by $Re = U D / \nu$ where $\nu$ is the kinematic viscosity.
%
%
Consider a sampling of trained Reynolds number values ranging from $90$ to $300$ with a jump equal to $30$. 
For each value, the numerical simulations were performed (Taylor-Hood finite element $\mbb{P}_2 / \mbb{P}_1 $) by using a time step $ dt = 0.01 $ and a non-uniform mesh including $ 85124$ DOFs for velocity and $ 10694$ DOFs for pressure. The initial condition considered for all the trained Reynolds number values is the solution at a given instant of the periodic flow regime at $ Re = 100 $. The final time of simulation was chosen equal to $12$.
\\
By considering $500$ snapshots regularly distributed in the time interval $[t_1, t_2] = [7,12]$, that is about $8$ periods of the flow, POD bases for velocity and pressure were calculated. The number of kept POD modes for velocity is equal to $10$ while $8$ modes are kept for pressure. 
In the following, we are going to demonstrate numerically on two optimization problems, the capability of RGA to reproduce in near real time a sufficiently accurate solution for the optimal control problem of flow past a cylinder when the Reynolds number value varies. For this optimal control problem, we mention that in the application of the Hyper Bi-CITSGM method inside the RGA, only three trained POD bases for which the corresponding parameters are the closest to the untrained parameter are considered in the ITSGM stage.   
Two minimization problems are separately studied in this example of the control of flow past a cylinder. The first problem is associated to the minimization of the following cost function 
\begin{equation}\label{EQ: Func_Strouhal_NS_GA}
 \functional^S(\statevecU,p) = 100 \times |S_t - \hat{S}_t | /  |\hat{S}_t |
\end{equation}
where $S_t$ is the Strouhal number. While the second problem is associated to the minimization of the cost function given by
\begin{equation}\label{EQ: Func_CLRMS_NS_GA}
\functional^L(\statevecU, p) = 100 \times | C_{L,rms} - \hat{C}_{L,rms} |/ |\hat{C}_{L,rms} | 
\end{equation}
where $C_{L,rms}$ is the root mean square lift coefficient defined by
\begin{equation*}
C_{L,rms} = \sqrt{\myfrac{1}{t_2-t_1} \int_{t_1}^{t_2} C_{L} \,dt}
\end{equation*}
The first (resp. second) control problem consists to determine by RGA the Reynolds number value that minimizes the above cost function \eqref{EQ: Func_Strouhal_NS_GA} (resp. \eqref{EQ: Func_CLRMS_NS_GA}).
%
%
In order to ensure that RGA delivers the best solution, we need to enrich the chromosomes by four additional genes. These are in this case the GIDW velocity and pressure weights powers $(l_u,m_u)$ and $(l_p,m_p)$ involved in equations \eqref{Bi-CITSGM_weights}. The space of search for RGA is given by
$$ K = \left\{(Re, l_u, m_u, l_p, m_p)\in \mbb{R}_+^5, \hspace*{0.3cm} 90 \leq Re \leq 300 \txt{ and } 1 < l_u, m_u, l_p, m_p \leq 8 \right\}$$
A population of $20$ chromosomes of $5$ genes randomly generated in $K$ is used as initial guess to obtain the numerical results; and our algorithm is run for $30$ generations. 
%
%
Consider the tests outlined in table \ref{tab:tests_cylinder}, where $\hat{S}_t$ (resp. $\hat{C}_{L,rms}$) is the target value associated to the minimization problem of the cost function \eqref{EQ: Func_Strouhal_NS_GA} (resp. \eqref{EQ: Func_CLRMS_NS_GA}) and $Re_{_{\txt{opt}}}$ is the corresponding Reynolds number value that has to be approximated by a value $Re_{_{\txt{GA}}}$ delivered at the end of RGA. 
\begin{table}[hbtp!]
\centering
 \begin{tabular}{ c|ccc }
   & $Re_{_{\txt{opt}}}$        & $\hat{S}_t$ &  $\hat{C}_{L,rms}$\\
 \hline
 Test 1	&	$135$ 	&	$0.165$	&	$0.347$		 \\
 Test 2	&	$160$	&	$0.170$	&	$0.408$		\\
 Test 3 &	$195$	&	$0.175$	&	$0.475$		 \\
 Test 4 &	$225$	&	$0.179$	&	$0.531$		 \\
\hline
\end{tabular}
\caption{Studied numerical tests for the control of flow past a cylinder using RGA. The values of $\hat{S}_t$ (resp. $\hat{C}_{L,rms}$) correspond to the target values considered in the minimization of $\functional^S$ (resp. $\functional^L$).}
\label{tab:tests_cylinder}
\end{table}
%
%
%
%
Figure \ref{average_func_cylinder} illustrates the decreasing behavior of the averaged cost functions $\bm{avg}(\functional^S)$ and $\bm{avg}(\functional^L)$. After $12$ generations, it can be observed that these averaged functions stagnate meaning by that the population contains a chromosome with high recurrence.
%
%
%
Table \ref{tab:Strouhal_and_CLRMS} presents the best chromosomes of generation $12$. 
It can be reported that RGA delivers a good approximation $Re_{_{\txt{GA}}}$ of the sought Reynolds number value $Re_{_{\txt{opt}}}$. 
Let's denote by $\bar{\varepsilon}^{\%}_f$ the percentage of error in the time interval $[t_1,t_2]$ between a given time dependent function $f$ and its approximation $\widetilde{f}$. $\bar{\varepsilon}^{\%}_f$ is defined by 
$$\bar{\varepsilon}^{\%}_f = 100 \times \left. \left( \int_{t_1}^{t_2}||f- \widetilde{f} ||_{L^2(\Omega)}^2  \, dt \right)^{\frac{1}{2}} \middle/ \left( \int_{t_1}^{t_2}|| f ||_{L^2(\Omega)}^2  \, dt \right)^{\frac{1}{2}} \right.$$
From table \ref{tab:Strouhal_and_CLRMS}, this percentage of error between the predicted RGA flow solution and the optimal solution is less than $2\%$ for velocity and $3.3\%$ for pressure. 
This can be further inspected visually from figure \ref{fig:3snapflow_RE160}, where a good match between the RGA and optimal solutions is observed.
In terms of CPU time, RGA reached a good approximation of the optimal Reynolds number value in about $32$ seconds, which proves the computational efficiency of this optimization approach.

\begin{figure}
\begin{subfigure}[b]{\textwidth}
\begin{tabular}{ c|cccccccc }
  &$Re_{_{\txt{opt}}}$& $Re_{_{\txt{GA}}}$  & $(l_u,m_u)$ & $(l_p,m_p)$     & $\functional^S(Re_{_{\txt{GA}}})$&$\bar{\varepsilon}^{\%}_u$&$\bar{\varepsilon}^{\%}_p$\\
 \hline
 Test 1	&$135$&	$135.4$ 	&$(4.20,3.05)$ & $(6.20,3.50)$ &	$0.10\%$&	$1.28 \%$ &	$ 3.10\%$ \\
 Test 2	&$160$&	$161.5$	& $(1.78,1.23)$ &$(2.82,5.28)$&	$0.12\%$&	$1.35 \%$ &	$2.79 \%$ \\ 
 Test 3 &$195$&	$193.1$	& $(4.74,1.70)$ &$(3.46,4.56)$& 	$0.001\%$&	$1.16 \%$ &	$ 3.30\%$ \\ 
 Test 4	&$225$&	$227.4$&	$(6.05,3.01)$	& $(4.55,3.40)$ & 	$0.0004\%$&	$1.56 \%$ &	 $3.35 \%$ \\ 
\hline
\end{tabular}
\caption{minimization of $\functional^S$}
\label{tab:Strouhal}
\end{subfigure}
\\
\vspace*{0.2cm}
\begin{subfigure}[b]{\textwidth}
\begin{tabular}{ c|cccccccc }
  &$Re_{_{\txt{opt}}}$& $Re_{_{\txt{GA}}}$  & $(l_u,m_u)$ & $(l_p,m_p)$     & $\functional^L(Re_{_{\txt{GA}}})$&$\bar{\varepsilon}^{\%}_u$&$\bar{\varepsilon}^{\%}_p$\\
 \hline
 Test 1	&$135$&	$136.8$ 	&$(4.60,4.75)$ & $(2.73,2.47)$ &	$0.25\%$&	$1.55\%$&	$ 2.94\%$ \\
 Test 2	&$160$&	$161.9$	& $(4.10,4.56)$ &$(5.92,2.74)$&	$0.01\%$&	$0.52 \%$ & $2.77 \%$ \\
 Test 3	&$195$&	$193.2$	& $(1.82,6.25)$ &$(5.97,5.20)$& 	$1.20\%$&	$0.74 \%$ & $3.13 \%$	 \\ 
 Test 4	&$225$&	$225.7$&	$(5.02,1.74)$	& $(2.35,2.02)$ & 	$0.17\%$&	$2.02 \%$ & $3.37 \%$	 \\
\hline
\end{tabular}
\caption{minimization of $\functional^L$}
\label{tab:CLRMS}
\end{subfigure}

\caption{Outputs at generation $12$ of RGA applied to the control problem of flow past a cylinder.}
\label{tab:Strouhal_and_CLRMS}
\end{figure}

%% file: Tex_Files/GA_Cavity.tex
\begin{figure}[hbtp!]
\begin{subfigure}{\textwidth}
\includegraphics[width=0.32\linewidth]{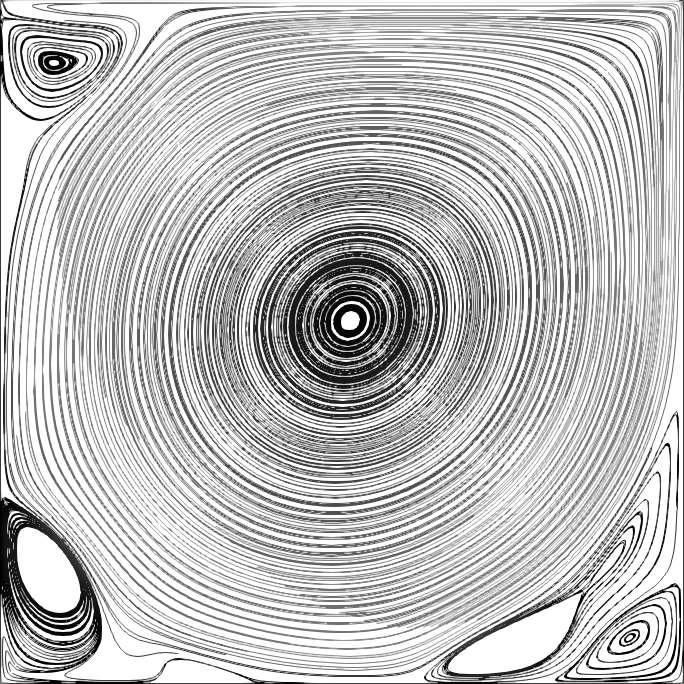}%
\hspace*{0.25cm}\includegraphics[width=0.32\linewidth]{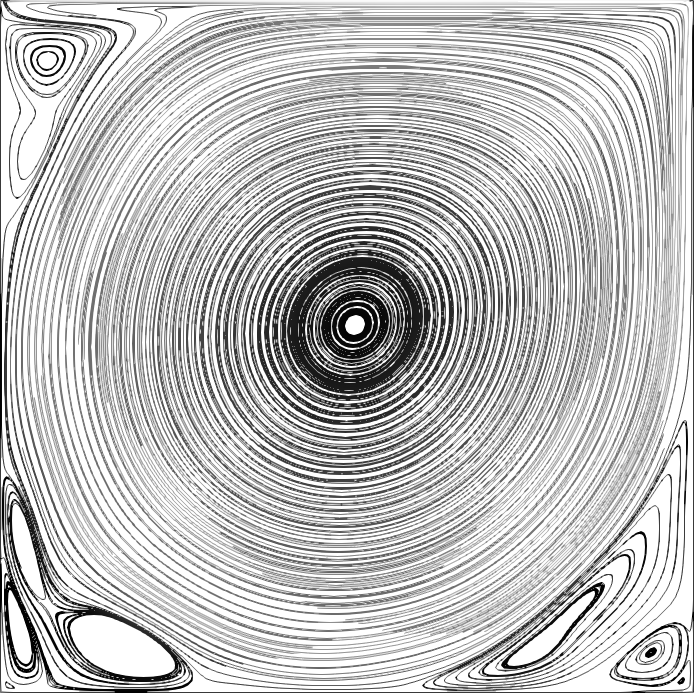}%
\hspace*{0.25cm}\includegraphics[width=0.32\linewidth]{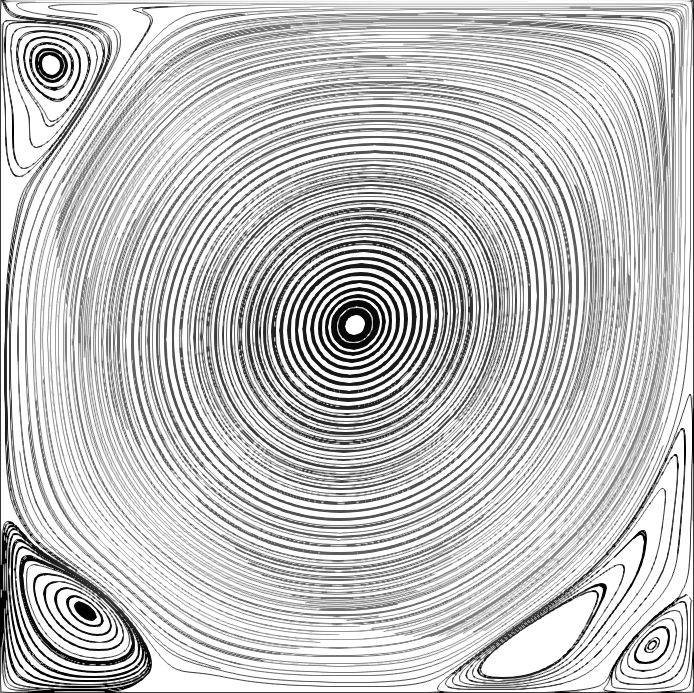}
\caption{flow solution at $Re = 9000$.}
\end{subfigure}
\begin{subfigure}{\textwidth}
\includegraphics[width=0.32\linewidth]{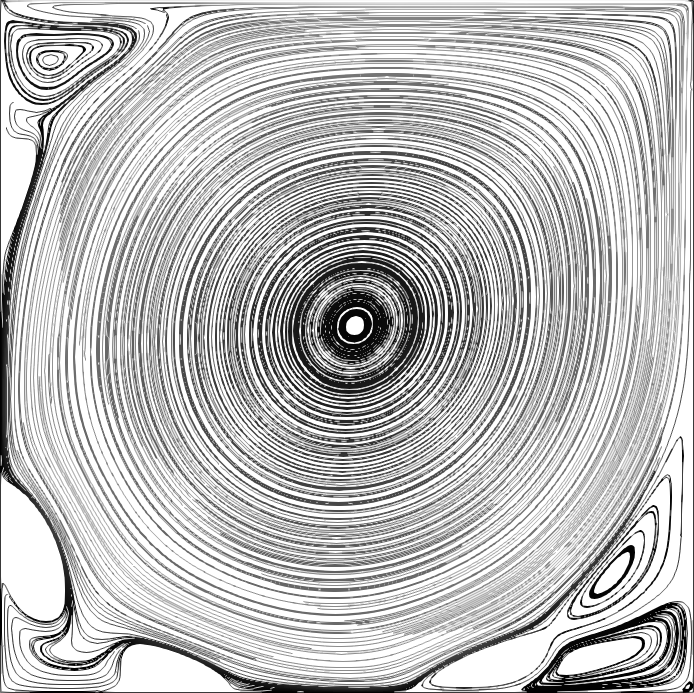}%
\hspace*{0.25cm}\includegraphics[width=0.32\linewidth]{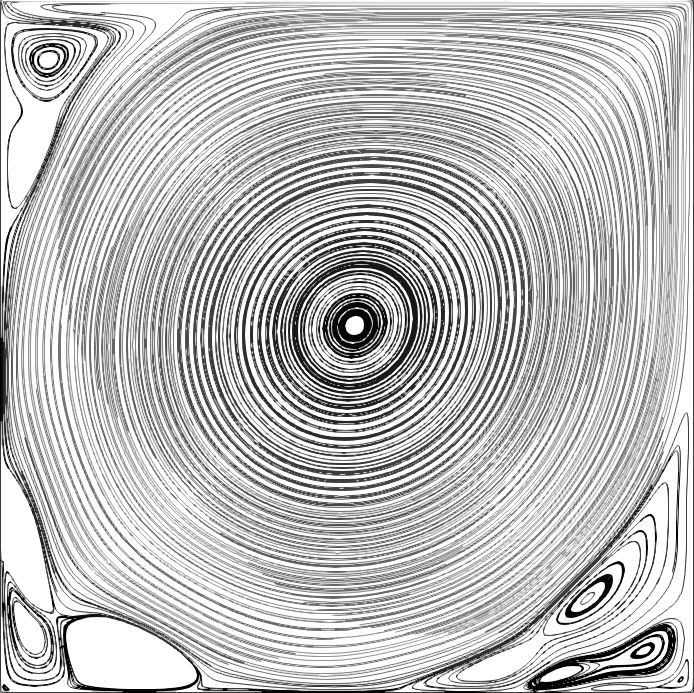}%
\hspace*{0.25cm}\includegraphics[width=0.32\linewidth]{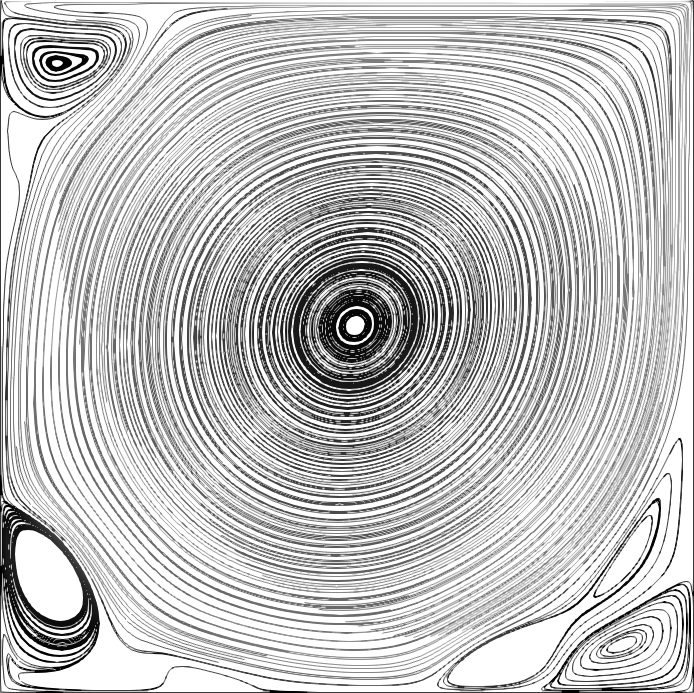}
\caption{flow solution at $Re = 10000$.}
\end{subfigure}
\begin{subfigure}{\textwidth}
\includegraphics[width=0.32\linewidth]{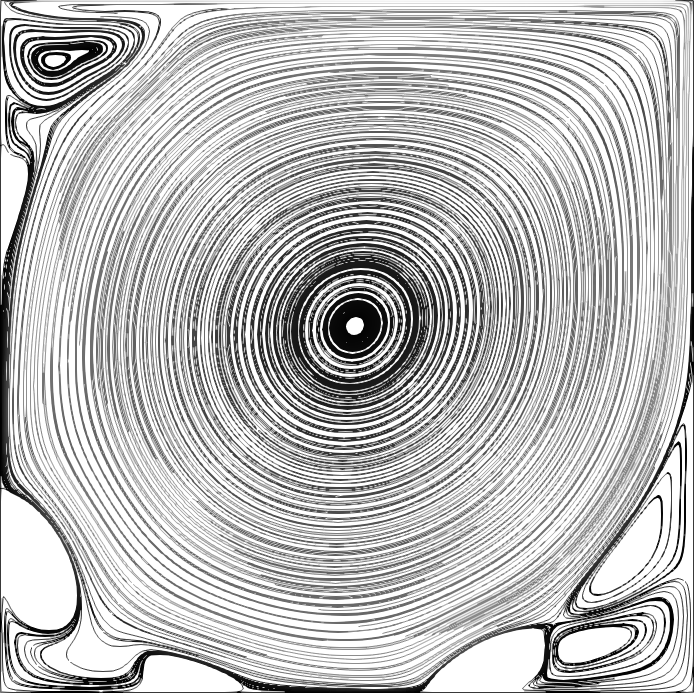}%
\hspace*{0.25cm}\includegraphics[width=0.32\linewidth]{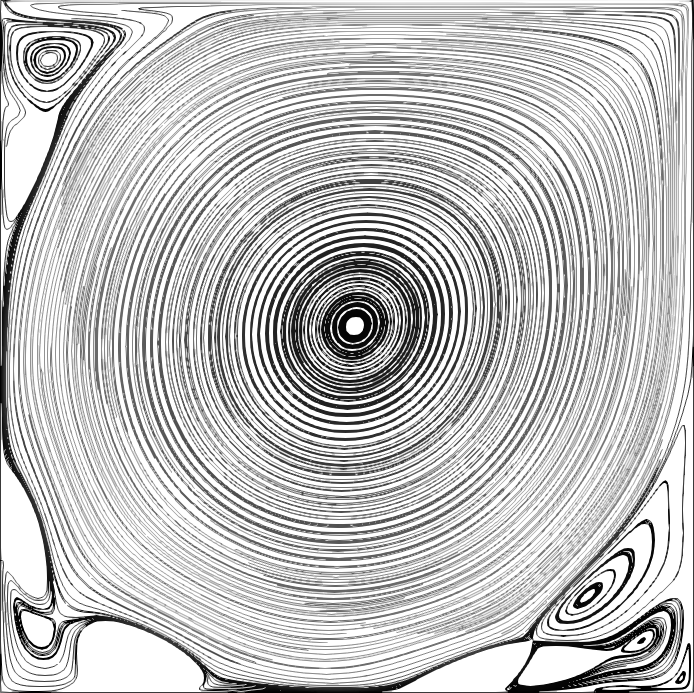}%
\hspace*{0.25cm}\includegraphics[width=0.32\linewidth]{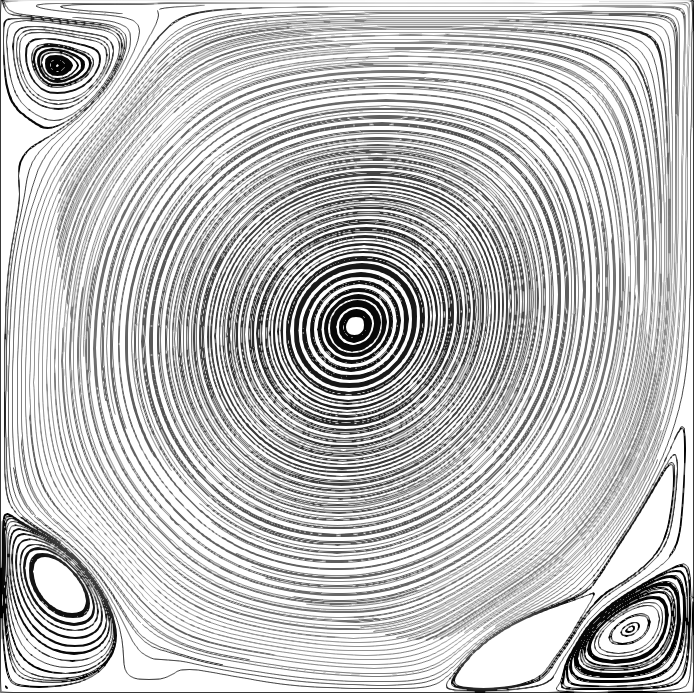}
\caption{flow solution at $Re = 11000$.}
\end{subfigure}
\caption{Evolution of the trained flow solutions at tree instants $t=T/4$ (left), $t=T/2$ (middle) and $t=T$ (right). Here, for each flow solution (at $Re = 9000$, $Re = 10000$ and $Re = 11000$), $T$ corresponds to the associated single flow period.}
\label{fig:Comparison_sampling_sol_cavity}
\end{figure}
\begin{figure}[hbtp!]
\hspace*{-2cm}

\centering
\includegraphics[width=0.6\linewidth]{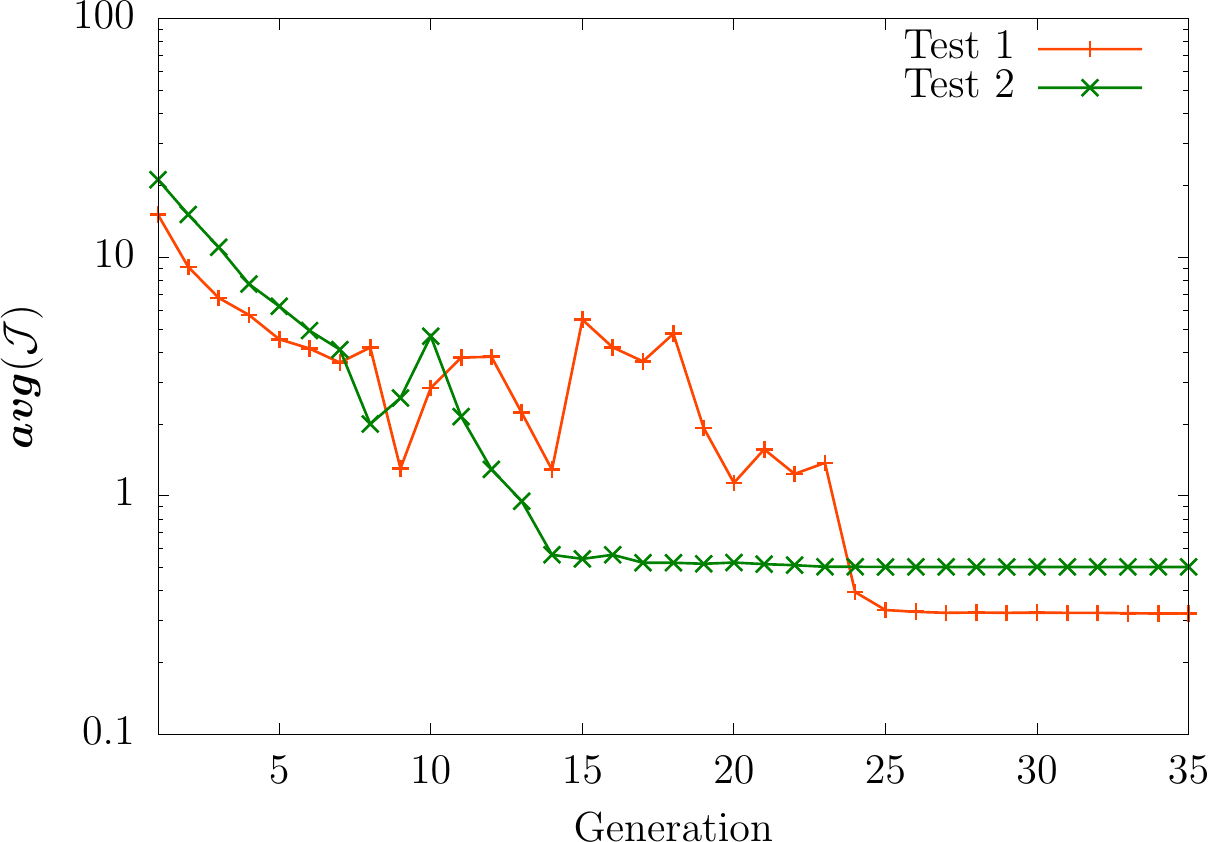}
\caption{Evolution of the averaged functional $\bm{avg}(\functional)$ with respect to generations. Test 1 corresponds to the optimal control $Re_{_{\txt{opt}}} = 9500$ and Test 2 to $Re_{_{\txt{opt}}} = 10500$.}
\label{average_func_cavity}
\end{figure}
\begin{figure}[hbtp!]
\begin{subfigure}{\textwidth}
\includegraphics[width=0.32\linewidth]{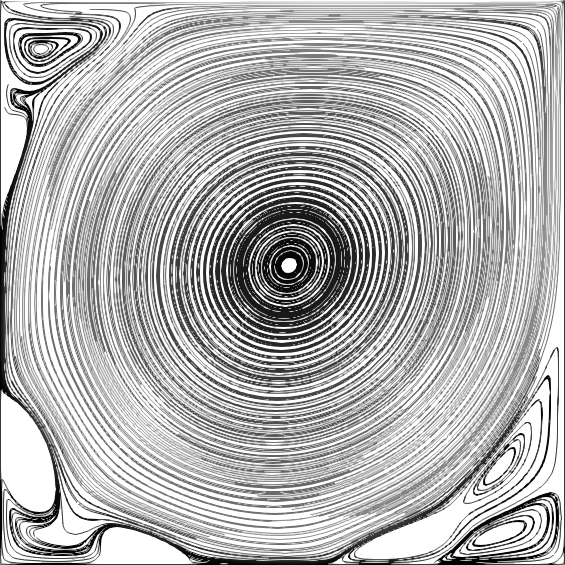}%
\hspace*{0.25cm}\includegraphics[width=0.32\linewidth]{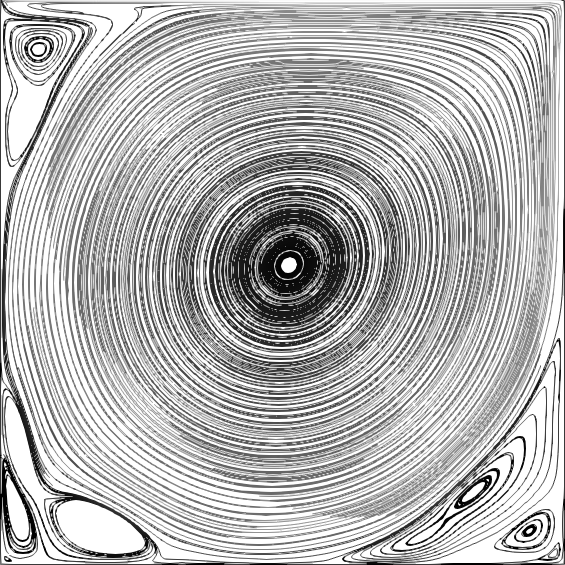}%
\hspace*{0.25cm}\includegraphics[width=0.32\linewidth]{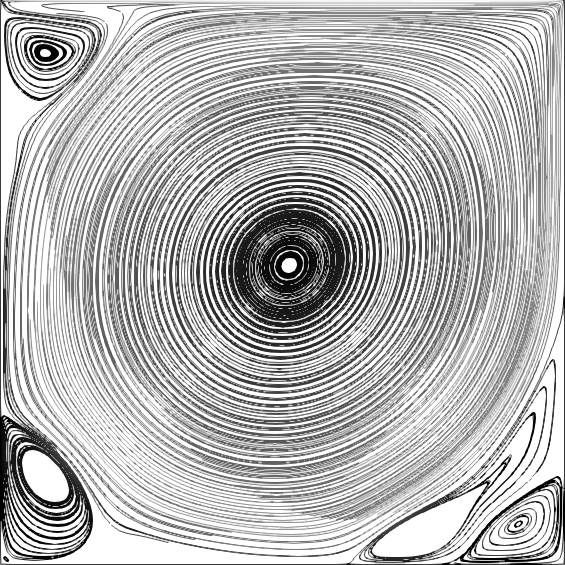}
\caption{Optimal flow solution}
\end{subfigure}
\begin{subfigure}{\textwidth}
\includegraphics[width=0.32\linewidth]{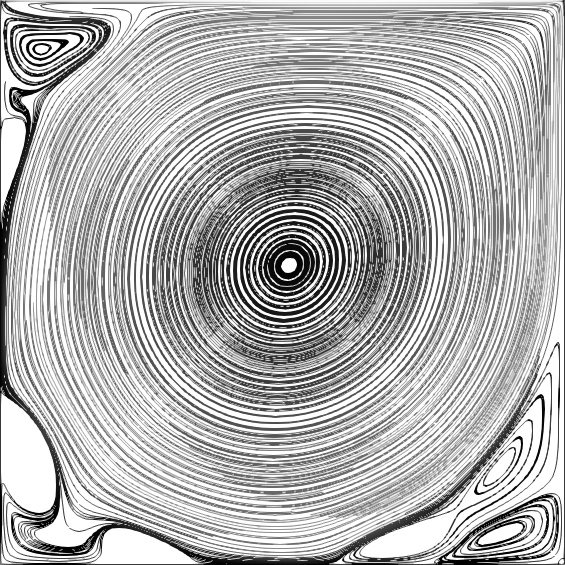}%
\hspace*{0.25cm}\includegraphics[width=0.32\linewidth]{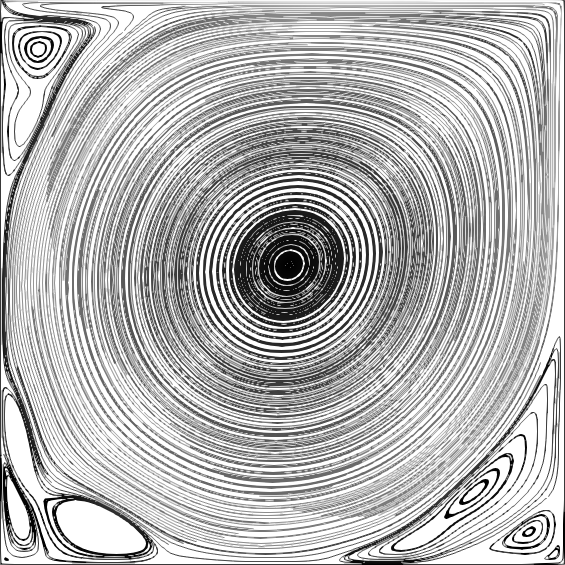}%
\hspace*{0.25cm}\includegraphics[width=0.32\linewidth]{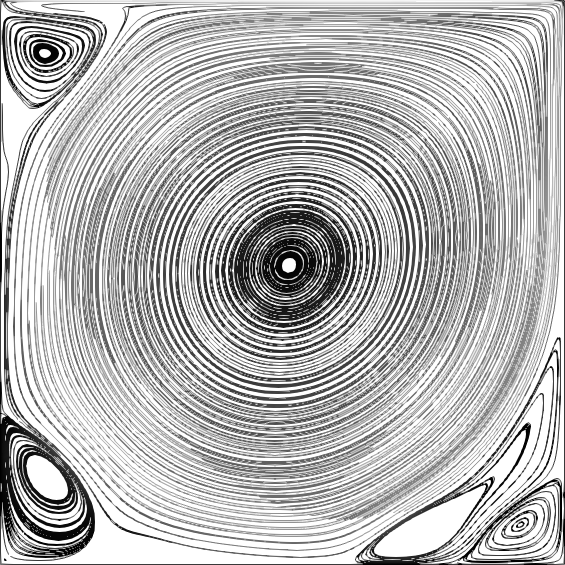}
\caption{flow solution obtained at the end of RGA applied to minimize $\functional$}
\end{subfigure}
\caption{Evolution of the target and predicted RGA solutions for Test 1 ($Re_{_{\txt{opt}}} = 9500$) at tree instants $t=T/4$ (left), $t=T/2$ (middle) and $t=T$ (right). Here, $T$ represents a single flow period at $Re = 9500$.}
\label{fig:Comparison_RGA_and_optimal_RE9500}
\end{figure}
\begin{figure}[hbtp!]
\begin{subfigure}{\textwidth}
\includegraphics[width=0.32\linewidth]{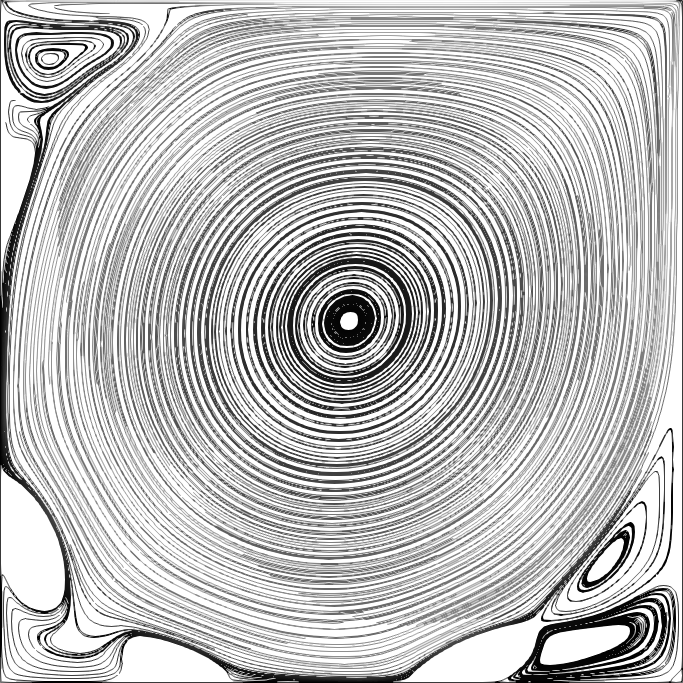}%
\hspace*{0.25cm}\includegraphics[width=0.32\linewidth]{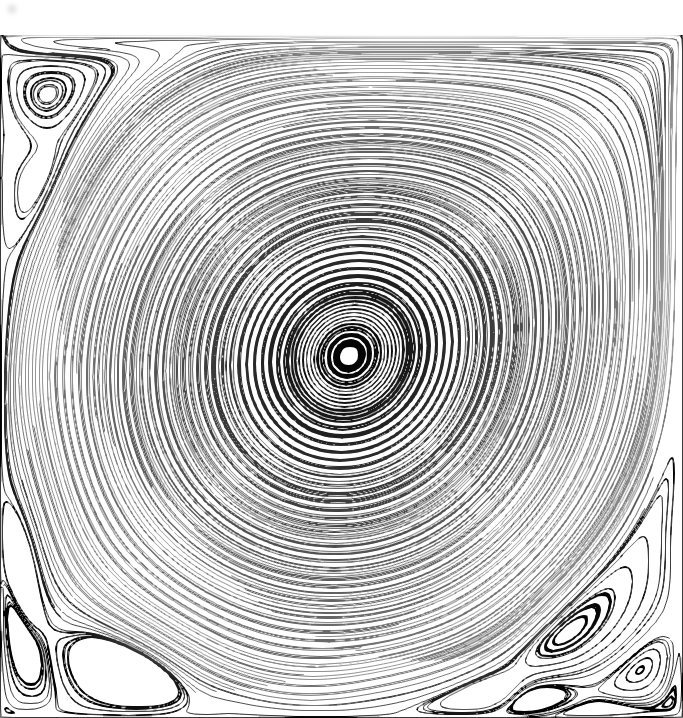}%
\hspace*{0.25cm}\includegraphics[width=0.32\linewidth]{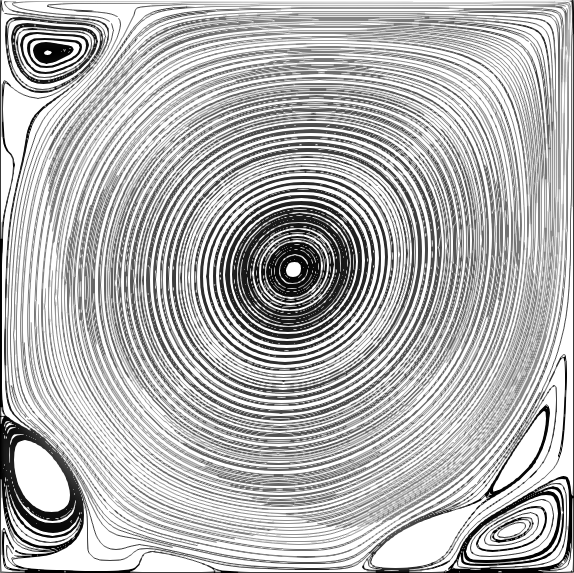}
\caption{Optimal flow solution}
\end{subfigure}
\begin{subfigure}{\textwidth}
\includegraphics[width=0.32\linewidth]{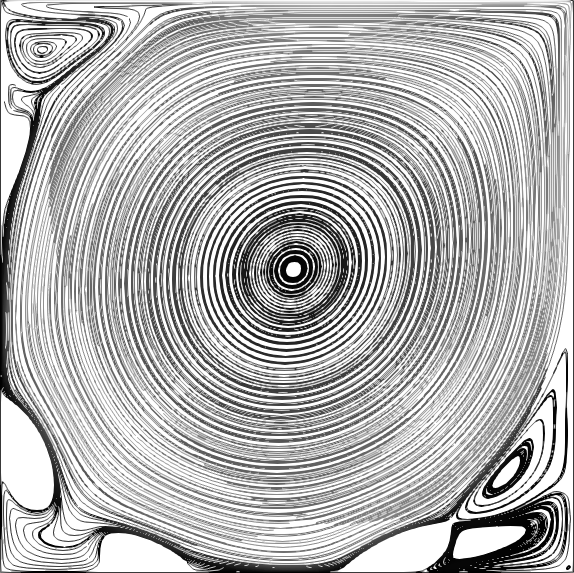}%
\hspace*{0.25cm}\includegraphics[width=0.32\linewidth]{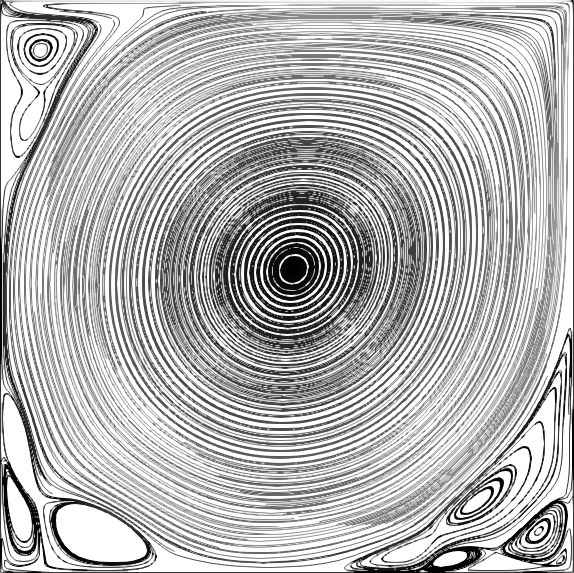}%
\hspace*{0.25cm}\includegraphics[width=0.32\linewidth]{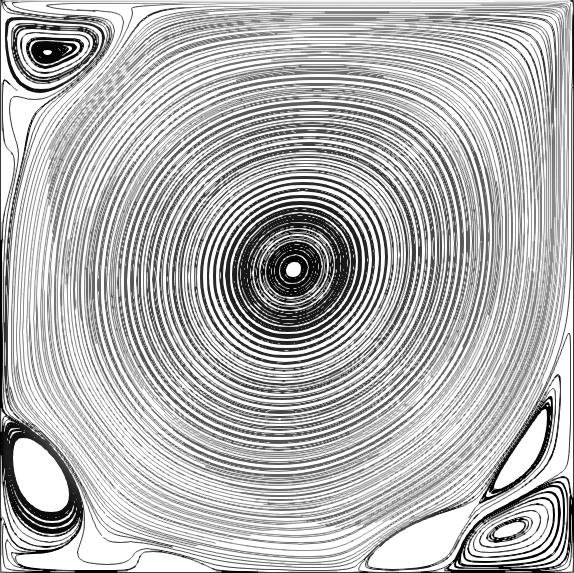}
\caption{flow solution obtained at the end of RGA applied to minimize $\functional$}
\end{subfigure}
\caption{Evolution of the target and predicted RGA solutions for Test 2 ($Re_{_{\txt{opt}}} = 10500$) at tree instants $t=T/4$ (left), $t=T/2$ (middle) and $t=T$ (right). Here, $T$ represents a single flow period at $Re = 10500$.}
\label{fig:Comparison_RGA_and_optimal_RE10500}
\end{figure}
\begin{figure}[hbtp!]
\hspace*{-2cm}
\begin{subfigure}{0.6\textwidth}
\includegraphics[width=\linewidth]{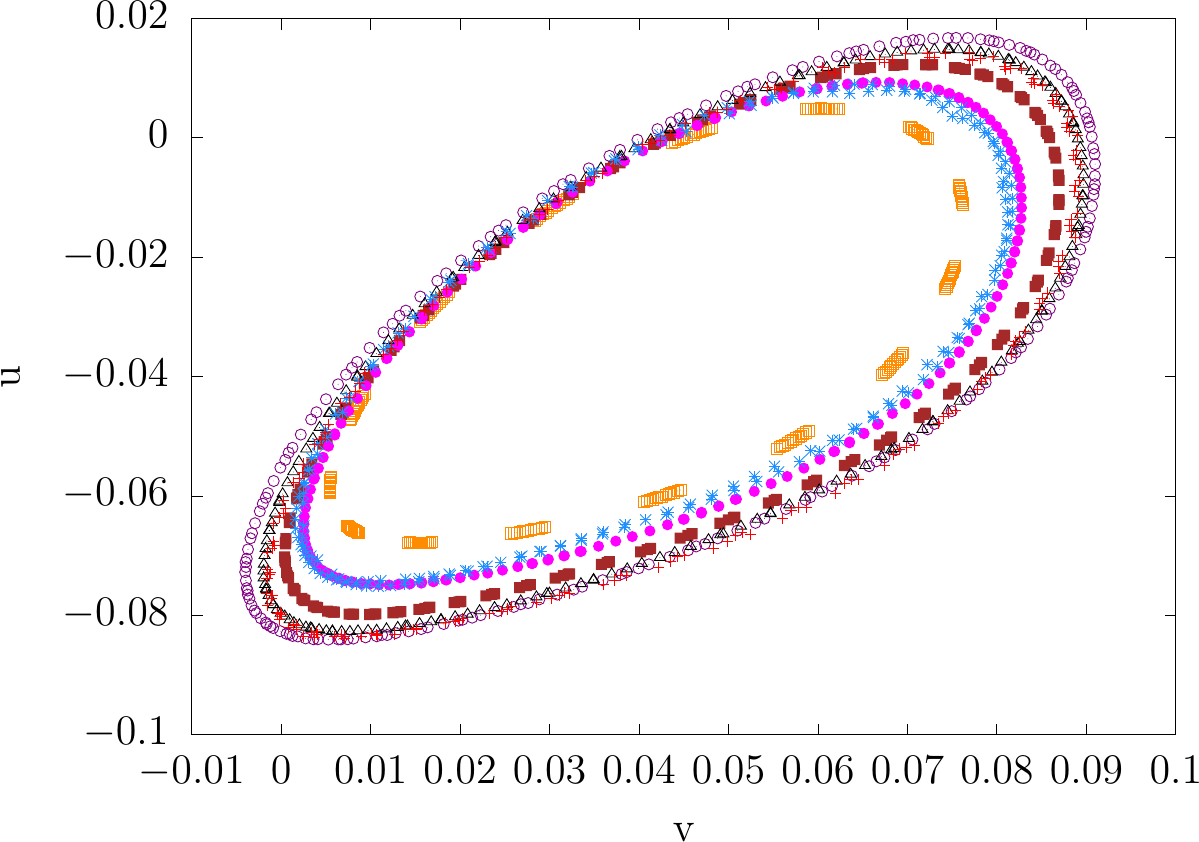}
\caption{Point $x_1 = (2/16, 13/16)$.}
\end{subfigure}%
\begin{subfigure}{0.6\textwidth}
\includegraphics[width=\linewidth]{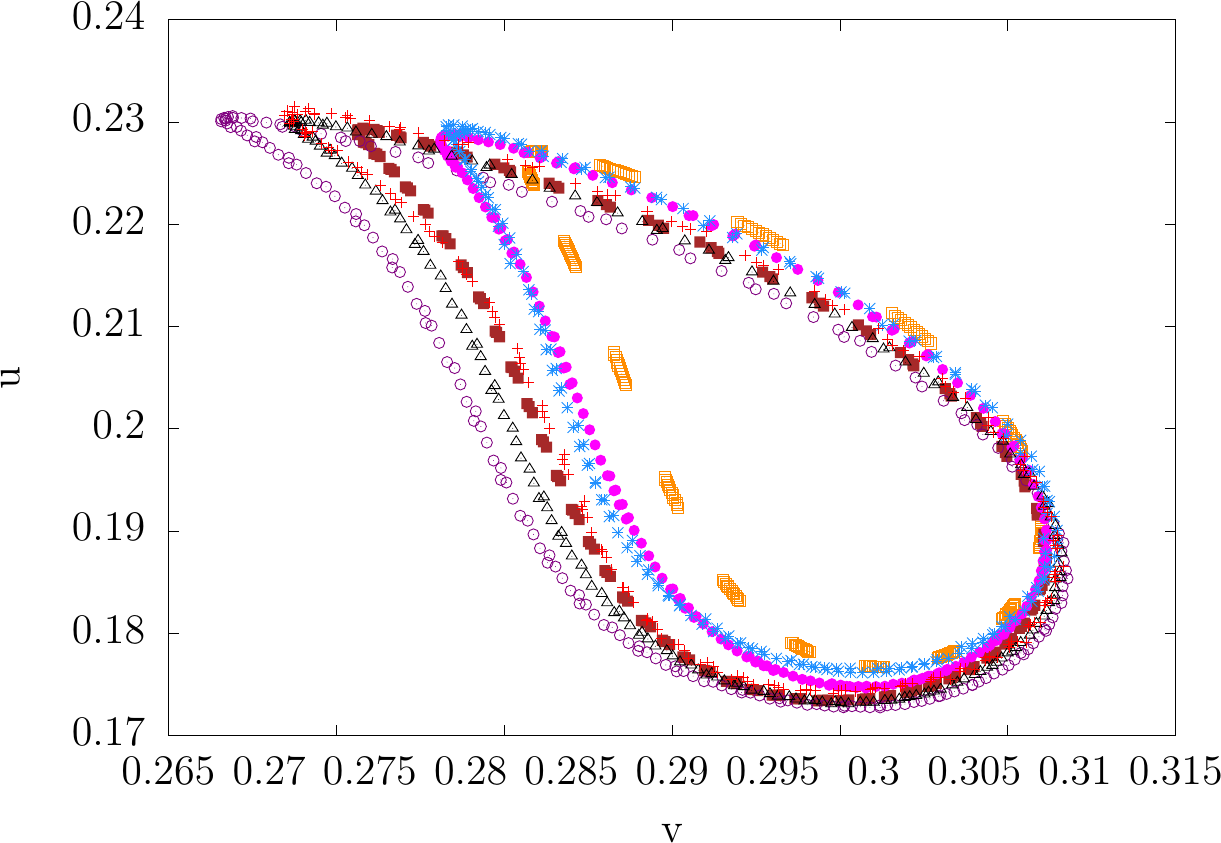}
\caption{Point $x_2 = (2/16, 2/16)$.}
\end{subfigure}
\hspace*{-2.0cm}
\begin{subfigure}{0.6\textwidth}
\includegraphics[width=\linewidth]{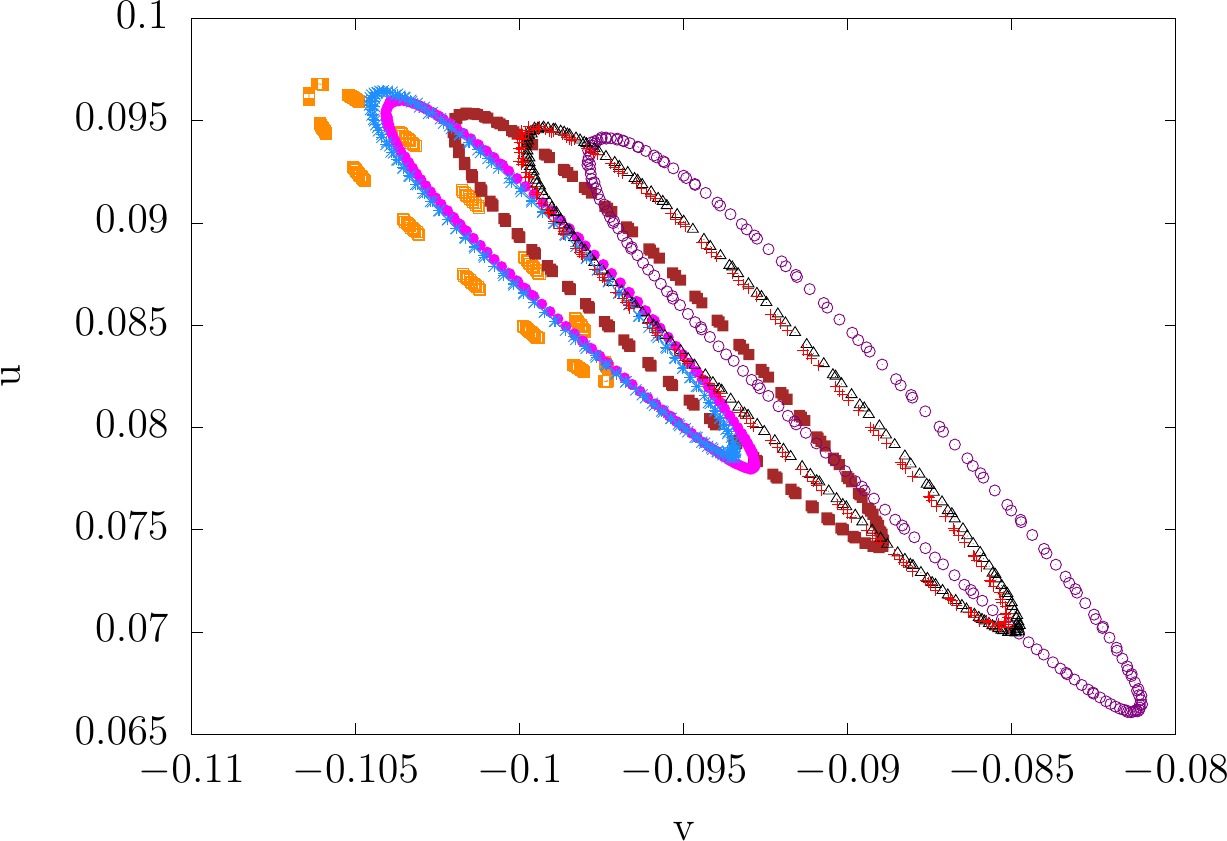}
\caption{Point $x_3 = (19/20, 19/20)$.}
\end{subfigure}%
\hspace*{0.2\textwidth}\begin{subfigure}{0.4\textwidth}
\includegraphics[width=\linewidth]{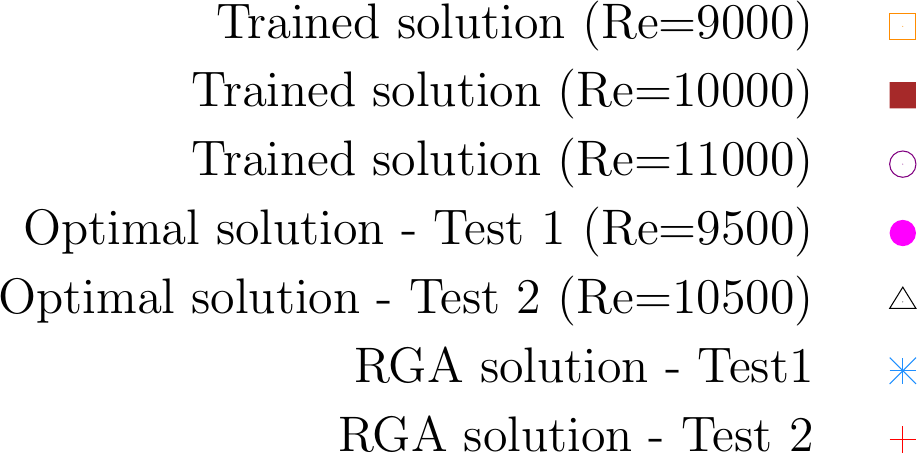}
\end{subfigure}
\caption{Velocity phase portrait at reference points $x_1$, $x_2$ and $x_3$ of the trained flow solutions (associated to $Re = 9000$, $Re=10000$ and $Re11000$) optimal solutions (associated to $Re = 9500$ and $Re=10500$) and RGA predicted solutions for Test 1 and Test 2. Here, $u$ and $v$ are respectively the horizontal and vertical components of the flow velocity $\bm{u}$.}
\label{fig:velocity_at_points}
\end{figure}

The two dimensional lid driven cavity flow problem is chosen now as test case for RGA. The problem domain consists of a square cavity $]0,D[\times]0,D[$ filled with fluid. At the top boundary, a tangential velocity $U$ of unit magnitude is applied to drive the fluid flow in the cavity, while the remaining three walls are defined as no-slip conditions. The Reynolds number of this flow is given by $Re = U D/\nu$.
Consider the trained Reynolds number values $9000$, $10000$ and $11000$. For each value, the nondimensional finite element high fidelity solver (Taylor-Hood finite element $\mbb{P}_2 / \mbb{P}_1 $) was performed on an unstructured mesh including $410370$ DOFs for velocity and $51585$ for pressure. The dynamics of the flow solutions at the chosen trained Reynolds number values can be observed from figure \ref{fig:Comparison_sampling_sol_cavity}, where snapshots at three different instants are represented. Based on the precalculated untrained flow solutions, $200$ snapshots uniformly selected from the periodic regime in the non dimensional time interval $[t_1,t_2]$ of length $20$, were used to build the velocity POD basis of dimension $8$. In the following, the goal is to apply the RGA to control the flow in a lid driven cavity by acting on the Reynolds number value. 
Consider the misfit function $\functional$ defined by the percentage error between the observed and calculated velocities at three points $x_1$,$x_2$ and $x_3$, i.e., 
\begin{equation*}
 \functional(\statevecU) = 100  \times \myfrac{\left(\somme{i}{1}{3} \int_{t_1+\tau}^{t_2-\tau}|\statevecU(t+\delta,x_i)-\statevecUCible(t,x_i)|^2 \,dt \right)^{1/2}}{\left(\somme{i}{1}{3} \int_{t_1+\tau}^{t_2-\tau} |\statevecUCible(t,x_i)|^2 \,dt\right)^{1/2}} 
\end{equation*}
where $\statevecUCible$ is the target velocity, $\delta$ is the phase shift such that $|\delta|\leq\tau$, and $x_1$,$x_2$ and $x_3$ are the nondimensional control points picked at the corners of the cavity where fluid recirculations are observed. These points are
\begin{equation*}
x_1 = (2/16, 13/16) \hspace*{1cm} x_2 = (2/16, 2/16) \hspace*{1cm} x_3 = (19/20, 19/20) \hspace*{1cm}
\end{equation*}
Chromosomes in RGA were enriched in this case by three additional genes. These are the velocity GIDW powers $(l_u,m_u)$ and the phase shift $\delta$.
The space of search is given as follows
$$ K = \left\{(Re, l_u, m_u, \delta)\in \mbb{R}_+^3\times\mbb{R}, \hspace*{0.3cm} 9000 \leq Re \leq 11000; 1 < l_u, m_u \leq 8 \hspace*{0.2cm}\txt{and}\hspace*{0.1cm} |\delta|\leq\tau \txt{ where } \tau =1 \right\}$$
A population of $70$ chromosomes of $4$ genes randomly generated in $K$ is used as initial guess to obtain the numerical results; and our algorithm is run for $100$ generations. 
Two numerical tests were considered by choosing target velocities $\statevecUCible$ associated respectively to $Re_{_{\txt{opt}}} = 9500$ (Test 1) and $Re_{_{\txt{opt}}} = 10500$ (Test 2).
Figure \ref{average_func_cavity} illustrates the decreasing behavior of the averaged cost function $\bm{avg}(\functional)$. It can be seen that the averaged RGA functional for Test 1, needed about $25$ generations to start stagnation, while for Test 2, the stagnation started earlier and $17$ generations were sufficient to declare a optimal control solution.
Table \ref{tab:Tracking_cavity} reports the best chromosomes as well as the consumed CPU time for the studied test cases. 
It can be reported from these results that RGA succeeded to delivers a good approximation $Re_{_{\txt{GA}}}$ of the sought Reynolds number value $Re_{_{\txt{opt}}}$. Moreover, the corresponding velocity solutions represented in figures \ref{fig:Comparison_RGA_and_optimal_RE9500} and \ref{fig:Comparison_RGA_and_optimal_RE10500} as well as the velocity phase portrait represented in figures \ref{fig:velocity_at_points}, show that RGA allowed a good prediction of the optimal flow solution for each test case. The percentage of error between the predicted velocity by RGA and the optimal solution was less than $1\%$. Finally, in terms of CPU time, we mention that RGA needed less than $16$ seconds to reach a good approximation of the sought optimal control, which proves once again the computational efficiency of the proposed optimization approach RGA.

\begin{table}[hbtp!]
 \centering
\begin{tabular}{ c|cccccccc }
 &$Re_{_{\txt{opt}}}$ &Generation& $Re_{_{\txt{GA}}}$  & $(l_u,m_u)$   & $\delta$     & $\functional(Re_{_{\txt{GA}}})$&$\bar{\varepsilon}^{\%}_u$&CPU time\\
 \hline
 Test 1	& $9500$ &$25$&	$9498.31$ 	&$(5.04, 4.18)$& $0.7$ &	$0.94\%$&$0.50\%$&	$16$ sec	 \\
 Test 2 	& $10500$&$17$&	$10502.76$	& $(2.62, 5.29)$& $0$ &	$1.34\%$&$0.72\%$&	$11$ sec	 \\
\hline
\end{tabular}
\caption{ Outputs and CPU time of RGA applied to the control problem of flow in a lid driven cavity.}
\label{tab:Tracking_cavity}
\end{table}
%
%
%
%

%% file: Tex_Files/Conclusion.tex
This paper presents a new application of the GA technique for near real time flow control. The proposed GA algorithm named RGA uses the hyper reduced version of Bi-CITSGM introduced in this study, which is a non intrusive model reduction approach that enables real time prediction of the flow solution.  
An attractive feature of RGA is its ability to evaluates and select chromosomes only based on a set of available trained high fidelity flow solutions data. In other words, no prior knowledge of the high fidelity mathematical equations is required. This implies that RGA can also be adopted to flow control for which data were obtained from experiment.
In this paper, the effectiveness of RGA was numerically tested on the control problems of flow past a cylinder and flow in a lid driven cavity where the control parameter was the Reynolds number value.
Interestingly, using a bunch of numerically precalculated trained high fidelity flow solutions, the RGA succeeded to provide good approximations of the sought optimal control solutions in less than half a minute. 
This shows that the proposed strategy is robust in terms of precision and CPU time for optimal control problems. 
Finally, it's worth noting that the time of RGA can be further reduced given the possibility of simultaneous parallel estimation of the cost functions for different control parameter combinations. Using such parallel implementation will speed up the optimal search by RGA and allow to achieve real-time optimal control solutions.

%% file: ARTICLE.bbl
\begin{thebibliography}{10}

\bibitem{OulghelouBiCITSGM2019Arxiv}
M.~Oulghelou and C.~Allery, ``Non intrusive method for parametric model order
  reduction using a bi-calibrated interpolation on the grassmann manifold,''
  {\em arXiv:1901.03177}, 2018.

\bibitem{Gunzburger2000}
M.~Gunzburger, ``Adjoint equation-based methods for control problems in
  incompressible, viscous flows,'' {\em Flow, Turbulence and Combustion},
  vol.~65, 12 2000.

\bibitem{Snyman2005}
J.~A. Snyman and D.~N. Wilke, {\em Practical mathematical optimization : an
  introduction to basic optimization theory and classical and new
  gradient-based algorithms}.
\newblock New York : Springer, 2005.

\bibitem{Desai1994}
M.~Desai and K.~Ito, ``Optimal controls of navier–stokes equations,'' {\em
  SIAM Journal on Control and Optimization}, vol.~32, no.~5, pp.~1428--1446,
  1994.

\bibitem{Zingg2008}
D.~W. Zingg, M.~Nemec, and T.~H. Pulliam, ``A comparative evaluation of genetic
  and gradient-based algorithms applied to aerodynamic optimization,'' {\em
  European Journal of Computational Mechanics}, vol.~17, no.~1-2, pp.~103--126,
  2008.

\bibitem{Holland1975}
J.~H. Holland, {\em Adaptation in Natural and Artificial Systems}.
\newblock Ann Arbor, MI: University of Michigan Press, 1975.
\newblock second edition, 1992.

\bibitem{Goldberg1989}
D.~E. Goldberg, {\em Genetic Algorithms in Search, Optimization and Machine
  Learning}.
\newblock Boston, MA, USA: Addison-Wesley Longman Publishing Co., Inc.,
  1st~ed., 1989.

\bibitem{SENGUPTA2007}
T.~K. Sengupta, K.~Deb, and S.~B. Talla, ``Control of flow using genetic
  algorithm for a circular cylinder executing rotary oscillation,'' {\em
  Computers \& Fluids}, vol.~36, no.~3, pp.~578 -- 600, 2007.

\bibitem{Hacioglu2005}
A.~Hac{\i}o\u{g}lu and \.{I}brahim \"{O}zkol, ``Inverse airfoil design by using
  an accelerated genetic algorithm via distribution strategies,'' {\em Inverse
  Problems in Science and Engineering}, vol.~13, no.~6, pp.~563--579, 2005.

\bibitem{SHAHROKHI2007}
A.~Shahrokhi and A.~Jahangirian, ``Airfoil shape parameterization for optimum
  navier–stokes design with genetic algorithm,'' {\em Aerospace Science and
  Technology}, vol.~11, no.~6, pp.~443 -- 450, 2007.

\bibitem{Daroczy2018}
L.~Dar\`{o}czy, G.~Janiga, and D.~Th\'{e}venin, ``Computational fluid dynamics
  based shape optimization of airfoil geometry for an h-rotor using a genetic
  algorithm,'' {\em Engineering Optimization}, vol.~50, no.~9, pp.~1483--1499,
  2018.

\bibitem{Vicini1999}
A.~Vicini and D.~Quagliarella, ``Airfoil and wing design through hybrid
  optimization strategies,'' {\em AIAA Journal}, vol.~37, no.~5, pp.~634--641,
  1999.

\bibitem{TerryAerodynamicSO}
Terry, L.~Hoist, Thomas, and H.~R. Pulliam, ``Aerodynamic shape optimization
  using a real-number-encoded genetic algorithm,''

\bibitem{Shigeru1997}
S.~Obayashi and T.~Tsukahara, ``Comparison of optimization algorithms for
  aerodynamic shape design,'' {\em AIAA Journal}, vol.~35, no.~8,
  pp.~1413--1415, 1997.

\bibitem{Hacioglu2002}
A.~Hac{\i}o\u{g}lu and \.{I}brahim \"{O}zkol, ``Vibrational genetic algorithm
  as a new concept in airfoil design,'' {\em Aircraft Engineering and Aerospace
  Technology}, vol.~74, no.~3, pp.~228--236, 2002.

\bibitem{Hacioglu2003}
A.~Hac{\i}o\u{g}lu and \.{I}brahim \"{O}zkol, ``Transonic airfoil design and
  optimisation by using vibrational genetic algorithm,'' {\em Aircraft
  Engineering and Aerospace Technology}, vol.~75, no.~4, pp.~350--357, 2003.

\bibitem{Doorly1998}
D.~J. Doorly and J.~Peir{\'o}, ``Supervised parallel genetic algorithms in
  aerodynamic optimisation,'' in {\em Artificial Neural Nets and Genetic
  Algorithms}, (Vienna), pp.~229--233, Springer Vienna, 1998.

\bibitem{Jones2000}
B.~R.~Jones, W.~A.~Crossley, and A.~Lyrintzis, ``Aerodynamic and aeroacoustic
  optimization of airfoils via a parallel genetic algorithm,'' 10 2000.

\bibitem{lumley1967}
J.~L. Lumley, ``{The Structure of Inhomogeneous Turbulent Flows},'' in {\em
  Atmospheric turbulence and radio propagation} (A.~M. Yaglom and V.~I.
  Tatarski, eds.), pp.~166--178, Moscow: Nauka, 1967.

\bibitem{Amsallem}
D.~Amsallem and C.~Farhat, ``An interpolation method for adapting reduced-order
  models and application to aeroelasticity,'' {\em AIAA Journal},
  pp.~1803--1813, 2008.

\bibitem{Lieu2004}
T.~Lieu and M.~Lesoinne, ``Parameter adaptation of reduced order models for
  three-dimensional flutter analysis,''

\bibitem{LIEU20065730}
``Reduced-order fluid/structure modeling of a complete aircraft
  configuration,'' {\em Computer Methods in Applied Mechanics and Engineering},
  vol.~195, no.~41, pp.~5730 -- 5742, 2006.
\newblock John H. Argyris Memorial Issue. Part II.

\bibitem{Bjorck71numericalmethods}
A.~Bj\:{o}rck and G.~H. Golub, ``Numerical methods for computing angles between
  linear subspaces,'' 1971.

\bibitem{OULGHELOUAMC2018}
M.~Oulghelou and C.~Allery, ``A fast and robust sub-optimal control approach
  using reduced order model adaptation techniques,'' {\em Applied Mathematics
  and Computation}, vol.~333, pp.~416 -- 434, 2018.

\bibitem{Oulghelou2017}
M.~Oulghelou and C.~Allery, ``Aip conference proceedings [author(s) icnpaa 2016
  world congress: 11th international conference on mathematical problems in
  engineering, aerospace and sciences - la rochelle, france (4–8 july 2016)]
  - optimal control based on adaptive model reduction approach to control
  transfer phenomena,'' vol.~1798, 2017.

\bibitem{Edelman1998}
A.~Edelman, T.~A. Arias, and S.~T. Smith, ``The geometry of algorithms with
  orthogonality constraints,'' {\em SIAM Journal on Matrix Analysis and
  Applications}, vol.~20, pp.~303--353, 01 1998.

\bibitem{Boumal-2015}
N.~Boumal and P.-A. Absil, ``Low-rank matrix completion via preconditioned
  optimization on the {G}rassmann manifold,'' {\em Linear Algebra and its
  Applications}, vol.~475, pp.~200--239, 06 2015.

\bibitem{Absil}
A.~Absil, R.~Mahony, and R.~Sepulchre, ``Riemann geometry of {G}rassmann
  manifolds with a view on algorithmic computation,'' {\em Acta Applicandae
  Mathematicae}, vol.~80, Issue 2, p.~199–220, 2004.

\bibitem{Wald}
R.~Wald, ``General relativity,'' {\em The University of Chicago Press}, 1984.

\bibitem{Sirovich}
L.~Sirovich, ``Turbulence and the dynamics of coherent structures : {P}art {I},
  {II} and {III},'' {\em Quarterly of Applied Mathematics}, pp.~461--590, 1987.

\bibitem{TheseRolando}
M.~M. Rolando, {\em Thèse : Interpolation sur les variétés {G}rassmanniennes
  et applications à la réduction de modèles en mécanique.}
\newblock Université de La Rochelle, 2018.

\bibitem{Rolando_IDW_article}
R.~Mosquera, A.~Hamdouni, A.~E. Hamidi, and C.~Allery, ``Pod basis
  interpolation via inverse distance weighting on grassmann manifolds,'' {\em
  Discrete \& Continuous Dynamical Systems - S}, vol.~0, p.~1743, 2018.

\end{thebibliography}
